\documentclass[11pt,reqno]{amsart}
\usepackage[T1]{fontenc}
\usepackage[utf8]{inputenc}
\usepackage{lmodern}
\usepackage[a4paper,margin=27mm]{geometry}
\usepackage{amsmath,amssymb,mathtools}
\usepackage[expansion=false]{microtype}
\usepackage{enumitem,needspace,graphicx,booktabs,tabularx,xurl}
\usepackage[hidelinks]{hyperref}
\numberwithin{equation}{section}
\newtheorem{theorem}{Theorem}[section]
\newtheorem{lemma}[theorem]{Lemma}
\newtheorem{proposition}[theorem]{Proposition}
\newtheorem{corollary}[theorem]{Corollary}
\theoremstyle{definition}
\newtheorem{definition}[theorem]{Definition}
\theoremstyle{remark}

\newtheorem{remark}[theorem]{Remark}

\newcommand{\R}{\mathbb R}
\newcommand{\C}{\mathbb C}
\newcommand{\Z}{\mathbb Z}
\newcommand{\N}{\mathbb N}
\newcommand{\Q}{\mathbb Q}
\newcommand{\T}{\mathbb T}
\newcommand{\I}{\mathbb I}
\newcommand{\Id}{\mathrm I}
\newcommand{\dd}{\,\mathrm d}

\newcommand{\norm}[1]{\left\lVert#1\right\rVert}

\newcommand{\distT}[1]{\left\lVert#1\right\rVert_{\T}}

\newcommand{\cZ}{\mathcal Z}

\DeclareMathOperator{\supp}{supp}
\DeclareMathOperator{\diag}{diag}
\DeclareMathOperator{\Leb}{Leb}

\setlist[enumerate]{label=\textup{(\arabic*)},leftmargin=*,itemsep=3pt}
\allowdisplaybreaks[2]
\title[Hausdorff dimension of the critical AMO spectrum]{Almost sure lower bound for the Hausdorff dimension of the spectrum of critical AMO}
\author{Yanhui Qu}
\address[Y.H.~Qu]{Department of Mathematical Sciences, Tsinghua University,\newline
  Beijing 100084, P.~R.~China.}
\email{yhqu@tsinghua.edu.cn}
\author{Pengyu Wang}
\address[P.~Y.~Wang]{Department of Mathematical Sciences, Tsinghua University,\newline
  Beijing 100084, P.~R.~China.}
\email{py-wang24@mails.tsinghua.edu.cn}
\date{}
\subjclass[2020]{Primary 47B36; Secondary 47A10, 28A78, 11K50}
\keywords{Critical almost Mathieu operator, Hausdorff dimension, spectral measure, continued fractions}
\hypersetup{bookmarksdepth=3,pdftitle={Almost sure lower bound for the Hausdorff dimension of the spectrum of critical AMO},pdfauthor={Yanhui Qu and Pengyu Wang},pdfkeywords={Critical almost Mathieu operator, Hausdorff dimension, spectral measure, continued fractions}}
\begin{document}
\begin{abstract}
We show that there exists an absolute constant $d>0$ such that for
almost all frequencies $\alpha\in[0,1]\setminus\Q$, the Hausdorff
dimension of the spectrum of the critical almost Mathieu operator
$H_{1,\alpha,\theta}$ is bounded from below by $d$.
\end{abstract}
\maketitle
\pagestyle{plain}
\section{Introduction}

\subsection{The context}

In this paper, we are interested in the Hausdorff dimension of the spectrum of the {\it almost Mathieu operator}  $H_{\lambda,\alpha,\theta}$ in $\ell^2(\mathbb Z)$  (denoted from now on   by AMO):
\begin{equation}\label{schro}
\ell^2(\mathbb Z) \ni u \mapsto (H_{\lambda,\alpha,\theta} u)_n= u_{n+1}+u_{n-1} +2\lambda \cos(
2\pi (n\alpha + \theta) )\, u_n\,,
\end{equation}
where $\theta\in \mathbb{T}=\R/\Z$ is  the phase, $\alpha\in (0,1)\backslash
\Q$ is the frequency and $\lambda>0$
is  the coupling constant.  
The AMO was first introduced by Peierls \cite{Pe},
as a model for an electron on a 2D lattice,  subjected to  a homogeneous
magnetic field \cite{Ha,R}. This model
has been extensively studied not only because of  its importance
in  physics \cite{AOS,OA,TKNN},
but also as a fascinating mathematical object.
For an account of the Hofstadter butterfly and its mathematical background,
see Bellissard's Bourbaki survey \cite{bel}.

The spectrum of
$H_{\lambda,\alpha,\theta}$ is a compact perfect set  in $\mathbb R$  and, observing that it is independent of $\theta$ since
$\alpha$ is irrational, we   denote it  by
$\Sigma_{\lambda,\alpha}$. For any $\lambda>0$, and for any irrational  $\alpha$,  $\Sigma_{\lambda,\alpha}$ is a Cantor set \cite{AJ05}. Each bounded connected component of  $\R\backslash\Sigma_{\lambda,\alpha}$ is called
an open gap   and we know from \cite{AJ08,AYZ2} that  if $\lambda \neq1$, then for any irrational   $\alpha$, all the spectral gaps allowed by gap labelling are open, as predicted by  the  Gap Labelling Theorem \cite{JM82}.

However, not much was known about the structure  of $\Sigma_{\lambda,\alpha}$ when $\lambda =1$.    $H_{1,\alpha,\theta}$ is called the {\it critical} AMO (or Harper's equation in physics);  it has particular importance in physics.
Helffer and Sj\"ostrand \cite{HS,HS3} developed a semiclassical analysis
of Harper's equation to study its spectral structure.  
{ One useful way to study $\Sigma_{\lambda,\alpha}$ is by periodic approximation \cite{AMS,k,last3},  through finer estimates of $\Sigma_{\lambda,p_n/q_n}$, where  $p_n(\alpha)/q_n(\alpha)$ is  the $n$th continued-fraction convergent of $\alpha$. Here $\Sigma_{\lambda,p/q}$ denotes the union of the spectra over all phases.
Helffer and Kerdelhu\'e \cite{HeKe} used semiclassical methods
to study the total bandwidth at rational frequencies. Last \cite{last3} shows that if $\alpha$ is not of bounded type, then $\Sigma_{1,\alpha}$ has  zero Lebesgue measure.  Finally, by the renormalization technique, Avila and Krikorian \cite{AK} completed the proof that $\Sigma_{1,\alpha}$ has zero Lebesgue measure, and hence is a Cantor set, for all irrational $\alpha$. Therefore, it is 
 natural to study the fractal dimensions of $\Sigma_{1,\alpha}$.}
  It was believed until the mid 1990's that
the box-counting  dimension  $\dim_B(\Sigma_{1,\alpha})$ equals $\frac{1}{2}$ for almost every $\alpha$; one can consult
\cite{conjhalf3,conjhalf2,conjhalf1}
for numerical and heuristic arguments supporting this conjecture. But  in 1994, Wilkinson-Austin \cite{wilkinson_austin}
provided  numerical evidence that $\dim_B(\Sigma_{1,\alpha})
\approx 0.498$\, for  $\alpha = \frac{\sqrt{5}-1}{2}$
and thus conjectured that $\dim_B(\Sigma_{1,\alpha}) < \frac{1}{2}$ for every irrational $\alpha$. However, Jitomirskaya-Zhang \cite{JZ} showed that  if $\beta(\alpha) > 0$, then $\overline{\dim}_B(\Sigma_{1,\alpha}) =1$, where $\overline{\dim}_B$ denotes upper box-counting dimension. This disproved Wilkinson-Austin's conjecture.
 Here, $\beta(\alpha)$,  which  measures how Liouvillean $\alpha$ is,  is defined as
\begin{equation}\label{defbeta}
\beta(\alpha):=\limsup_{n\rightarrow \infty}\frac{\log
q_{n+1}(\alpha)}{q_n(\alpha)},\end{equation}
where  $q_n(\alpha)$ is the denominator of  the     $n$th convergent of $\alpha.$
 Last  \cite{last3} showed that if  $ q_{n+1}(\alpha)>q_n^4(\alpha)$ for a subsequence  of $n$, then $\dim_H(\Sigma_{1,\alpha})\leq \frac{1}{2}$. We note that the set of such kind of frequencies is a
dense $G_\delta$ set which contains  $\{\alpha\in (0,1)\setminus\Q\, | \, \beta(\alpha)>0 \}$.   Last-Shamis \cite{LS}  showed that there exists a dense subset of $\{\alpha \in (0,1)\setminus\Q | \beta(\alpha)>0 \}$, for which $\dim_H(\Sigma_{1,\alpha})=0$. Recently,  Avila-Last-Shamis-Zhou \cite{ALSZ} strengthened the result of \cite{LS} and showed that for any $\alpha$  in   $ \{\alpha \in (0,1)\setminus\Q \, | \, \beta(\alpha)>0 \}$,  $\dim_H(\Sigma_{1,\alpha})=0$. As we can see, the  results of \cite{ALSZ,JZ,last3,LS} concern frequencies satisfying particular rational approximation conditions.
Recently, Becker, Jitomirskaya and Krasovsky \cite{BJK} made a breakthrough by showing that for any $\alpha\in[0,1]\setminus\Q$, the Hausdorff dimension of $\Sigma_{1,\alpha}$ is at most $1/2$.
The same authors subsequently showed that
$\overline{\dim}_B\Sigma_{1,\alpha}\le 1/2$ for Lebesgue almost every $\alpha$
\cite[Theorem~1.2]{BJKHarper}. 

The results above on Hausdorff dimension concern upper bounds. How about the lower bounds? In this direction, there is essentially only one existing result.  Helffer, Liu, Qu and Zhou \cite{HLQZ} showed that there exists a dense set of frequencies  $\mathcal F\subset [0,1]\setminus\Q$ (with positive Hausdorff dimension, yet  zero Lebesgue measure), such that $\dim_H \Sigma_{1,\alpha} > 0$ for $\alpha\in \mathcal F.$ Later, the same authors \cite{HLQZ25} strengthened the result by showing that for any $\delta>0$, there exists a set $\mathcal F_\delta\subset \mathcal F$ (with positive Hausdorff dimension) such that $0<\dim_H \Sigma_{1,\alpha}\le \overline{\dim}_B \Sigma_{1,\alpha}\le \delta$ for any $\alpha\in \mathcal F_\delta.$

\subsection{Main result}

In this paper, we obtain a uniform positive lower bound for the Hausdorff  dimension of $\Sigma_{1,\alpha}$ for a.e. $\alpha\in [0,1]\setminus \Q.$

\begin{theorem}\label{thm:main}
  There exists an absolute constant $d>0$ such that for Lebesgue a.e. $\alpha\in [0,1]\setminus \Q$,
  $$
  \dim_H \Sigma_{1,\alpha}\ge d.
  $$
\end{theorem}

\begin{remark}
  {\rm 
  (i)  Bellissard\footnote{ Private conversation with Y. Last, circ. 1995.} conjectured that there
exists some $\kappa \in (0, 1/2]$ such that
$\dim_H\Sigma_{1,\alpha} = \kappa$ for almost every
$\alpha$. Our result can be viewed as a small step towards this conjecture.

(ii) Based on the results concerning the Sturmian Hamiltonian in \cite{CQ2025}, we even expect that there
exists some $\kappa' \in (0, 1/2]$ such that
$\dim_H\mathcal N_{1,\alpha} = \kappa'$ for almost every
$\alpha$, where $\mathcal N_{1,\alpha}$ is the density of states measure of $H_{1,\alpha,\theta}$.

}
\end{remark}

Let us explain briefly the idea of the proof. By the chiral transformation introduced in \cite{BJK}, the critical AMO $H_{1,\alpha,\theta}$ has the same spectrum as a singular Jacobi operator $J_\alpha$. Let $J_\alpha^+$ be the half-line restriction of $J_\alpha$. It is not hard to observe that $\sigma(J_\alpha^+)=\sigma(J_\alpha)\mathrel{=}\Sigma_{1,\alpha}$. On the other hand, for almost every $\alpha$, we can show that the spectral measure $\mu_\alpha$ of $J_\alpha^+$ at $\delta_0$ is $d$-H\"older continuous. Since $\mu_\alpha$ is supported on $\sigma(J_\alpha^+)=\Sigma_{1,\alpha}$, by a well-known result from fractal geometry, we conclude that $\dim_H\Sigma_{1,\alpha}\ge d.$

From now on, we write 
$$ \I:=[0,1]\setminus \Q.$$
 Since in the rest of the paper we only deal with the critical AMO, for simplicity we write $$\Sigma_\alpha:=\Sigma_{1,\alpha}.$$

For nonnegative quantities, $f\lesssim_\alpha g$ means that $f\le C(\alpha)g$ with a constant depending only on $\alpha$. Other subscripts indicate the same convention, and $g\gtrsim_\alpha f$ is equivalent to $f\lesssim_\alpha g$.
\par\medskip
\noindent\textbf{Use of AI tools.}
The authors used ChatGPT-6 Astra and Codex (OpenAI) for language editing and to assist in exploring and developing
parts of the proof of Lemma~\ref{lem:recursive}, particularly the local estimate
in Lemma~\ref{prop:area}. The authors critically assessed the suggestions and
reconstructed and revised the proofs. They take full responsibility
for the mathematical arguments and the content of the final manuscript.

\section{Proof of Theorem~\ref{thm:main}}\label{sec:mainproof}
\subsection{Definition and basic properties of \texorpdfstring{$J_{\alpha,\theta}$}{the Jacobi operator}}
For any $\alpha\in\I$, $\theta\in\T$ {with a chosen real representative}, and $n\in\Z$, define
\begin{equation}\label{eq:bcoeff}
 b_n(\alpha;\theta):=2\sin\pi(n\alpha+\theta)
 \quad\text{and}\quad b_n(\alpha):=b_n(\alpha;0)=2\sin(n\pi\alpha).
\end{equation}
Define the Jacobi operator $J_{\alpha,\theta}:\ell^2(\Z)\to\ell^2(\Z)$ as
\begin{equation}\label{eq:Jacobi}
 (J_{\alpha,\theta}\phi)_n:=b_n(\alpha;\theta)\phi_{n-1}
      +b_{n+1}(\alpha;\theta)\phi_{n+1}\qquad(n\in\Z).
\end{equation}
When $\theta=0$, we write
\begin{equation}\label{eq:J}
 J_\alpha:=J_{\alpha,0}.
\end{equation}
Define $J_{\alpha,\theta}^+$ as the restriction of $J_{\alpha,\theta}$ to
$\ell^2(\Z^+)$ with Dirichlet boundary condition.
{ Here $\Z^+=\{0,1,\ldots\}$ and $J_\alpha^+:=J_{\alpha,0}^+$.
Changing the representative of $\theta$ by one changes the sign of every coefficient.
Multiplication of the $n$th coordinate by $(-1)^n$ intertwines the two operators.}

\subsection{Proof of Theorem~\ref{thm:main}}
For any $\alpha\in\I$, let $\mu_\alpha$ be the spectral measure of $J_\alpha^+$
{at the coordinate vector $\delta_0$}. That is, $\mu_\alpha$ is the unique
probability measure on $\sigma(J_\alpha^+)$ such that for any
{bounded Borel function $f$ on $\sigma(J_\alpha^+)$},
\begin{equation}\label{eq:mu}
 {\langle\delta_0,f(J_\alpha^+)\delta_0\rangle}
 =\int_{\sigma(J_\alpha^+)} f(E)\dd\mu_\alpha(E).
\end{equation}
{We take the inner product to be linear in its second argument.}
\begin{lemma}\label{lem:support}
For any $\alpha\in\I$, one has
\[
 \sigma(J_\alpha)=\sigma(J_\alpha^+)
 =\supp\mu_\alpha=\Sigma_\alpha.
\]
\end{lemma}
{
\begin{proof}
{
In the notation of \cite[equation~(2)]{BJK}, the coefficient joining $j$
and $j+1$ in $\widehat H_{\beta,\xi}$ is $2\sin(2\pi(\beta j+\xi))$.
For irrational $\alpha$, the chiral transform gives
$\sigma(\widehat H_{\alpha/2,\xi})=\Sigma_\alpha$ for every $\xi$, as stated
in the proof of \cite[Theorem~7.11]{BJK}. Taking $\xi=\alpha/2$ gives
the coefficient $b_{j+1}(\alpha)$, so that
\[
 J_\alpha=\widehat H_{\alpha/2,\alpha/2},
 \qquad \sigma(J_\alpha)=\Sigma_\alpha .
\]}

Since $b_0(\alpha)=0$, the two half-lines reduce $J_\alpha$. Reflection
$j=-1-r$ identifies the negative half with $-J_\alpha^+$, because
$b_{-r-1}=-b_{r+1}$. The alternating-sign unitary on $\ell^2(\Z^+)$
conjugates $J_\alpha^+$ to $-J_\alpha^+$. Hence
$\sigma(J_\alpha^+)=\sigma(J_\alpha)$. Cyclicity of the endpoint
vector $\delta_0$ for $J_\alpha^+$ now gives
$\supp\mu_\alpha=\sigma(J_\alpha^+)=\Sigma_\alpha$, as claimed.
\end{proof}}
\begin{lemma}\label{lem:holder}
There exists an absolute constant $d>0$ such that for Lebesgue a.e.\
$\alpha\in\I$ there exists a constant $C(\alpha)>0$ such that for any
$E\in\R$ and $r\in(0,1)$,
\[
 \mu_\alpha(B(E,r))\le C(\alpha)r^d .
\]
\end{lemma}

We will prove Lemma~\ref{lem:holder} in Section~\ref{sec:measure}.

\begin{proof}[Proof of Theorem~\ref{thm:main}]
Let $\widetilde{\I}\subset\I$ be the full measure set of frequencies in Lemma~\ref{lem:holder}.
Fix any $\alpha\in\widetilde{\I}$. On one hand, $\mu_\alpha$ is supported on
$\Sigma_\alpha$ by Lemma~\ref{lem:support}. On the other hand, for any $E\in\Sigma_\alpha$,
by Lemma~\ref{lem:holder}, the lower local dimension of $\mu_\alpha$ at $E$ satisfies
\[
 \underline d_{\mu_\alpha}(E)=\liminf_{r\downarrow0}
 \frac{\log\mu_\alpha(B(E,r))}{\log r}\ge d .
\]
Then by \cite[Proposition~2.3(a)]{Falconer1997}, $\dim_H\Sigma_\alpha\ge d$.
\end{proof}
\section{Proof of Lemma~\ref{lem:holder}}\label{sec:measure}
\subsection{Two special formal solutions}\label{subsec:solutions}\label{subsec:secondmatrix}
For any $\alpha\in\mathbb I$ and $z\in\C$, define two special formal solutions
$u(z,\alpha),v(z,\alpha)\in\C^{\Z^+}$ as follows:
\begin{equation}\label{eq:def-u-v}
 \begin{cases}
  (J_\alpha^+-z)u(z,\alpha)=0,&u_0(z,\alpha)=1,\\
  (J_\alpha^+-z)v(z,\alpha)=\delta_0,&v_0(z,\alpha)=0.
 \end{cases}
\end{equation}
These are formal sequences; they need not belong to $\ell^2$.
When $\alpha$ is fixed we omit it from the notation. More precisely,
\begin{align}
 &u_0=1,\quad u_1=z/b_1,\quad
 b_{j+1}u_{j+1}=zu_j-b_ju_{j-1}\quad(j\ge1),\label{eq:prec}\\
 &v_0=0,\quad v_1=1/b_1,\quad
 b_{j+1}v_{j+1}=zv_j-b_jv_{j-1}\quad(j\ge1).\label{eq:vrec}
\end{align}
For any $n\in\N$ and $\phi\in\C^{\Z^+}$, define
\[
 \phi^{[n]}:=\Pi_n(\phi)=(\phi_0,\ldots,\phi_{n-1})^T .
\]
{ For a finite index set $S$, $\Pi_S$ restricts to the coordinates in $S$,
listed in increasing order. In a product acting on the original coordinate space,
we identify the restricted vector with its zero extension. Thus $\Pi_n$
selects the initial segment $\{0,\ldots,n-1\}$. For an indexed matrix,
$A^{[n]}$ is its leading $n\times n$ block, and $\|\cdot\|$ denotes the
Euclidean norm or its induced operator norm.}
\begin{lemma}\label{lem:solutions}
Fix any $\alpha\in\mathbb I$ and $z\in\C$. Then:
\begin{enumerate}
\item For every $n\in\Z^+$,
$|u_n(z,\alpha)|\ge|u_n(\operatorname{Re}z,\alpha)|$.
\item There exists a lower triangular matrix $H(\alpha)$ such that,
for every positive integer $n$,
\[
 v^{[n]}(z,\alpha)=H^{[n]}(\alpha)u^{[n]}(z,\alpha).
\]
\item For any $n\ge1$,
\[
 W(u(z,\alpha),v(z,\alpha))(n):=
 u_n(z,\alpha)v_{n-1}(z,\alpha)-u_{n-1}(z,\alpha)v_n(z,\alpha)
 =-\frac{1}{b_n(\alpha)} .
\]
\end{enumerate}
\end{lemma}
\begin{proof}
(1) Let $\chi_j(z)=\det(z\Id_j-(J_\alpha^+)^{[j]})$ and $\chi_0=1$.
Expansion along the last row gives
$\chi_{j+1}=z\chi_j-b_j^2\chi_{j-1}$, with $\chi_1=z$.
Comparison with \eqref{eq:prec} yields
$u_j(z)=\chi_j(z)/(b_1\cdots b_j)$.
Since $(J_\alpha^+)^{[j]}$ is real symmetric, all its eigenvalues
$\xi_1,\ldots,\xi_j$ are real. Writing $z=E+\mathrm it$ gives
\[
 |u_j(E+\mathrm it)|^2
 =|b_1\cdots b_j|^{-2}
   \prod_{s=1}^j\bigl((E-\xi_s)^2+t^2\bigr)
 \ge|u_j(E)|^2.
\]
The case $j=0$ is immediate.

(2) For $i,j\ge0$, write
\begin{equation}\label{eq:matrix-indices}
 \kappa_{ij}:=\frac{i+j+1}{2},\qquad \sigma_{ij}:=\frac{i-j-1}{2}.
\end{equation}
Define the strictly lower triangular matrix
\begin{equation}\label{eq:Hentries}
 H_{ij}(\alpha)=
 \begin{cases}
 \displaystyle\frac{(-1)^{\sigma_{ij}}}{b_{\kappa_{ij}}},
       &i>j,\quad i-j\text{ odd},\\[5pt]
 0,&\text{otherwise}.
 \end{cases}
\end{equation}
Put $w_i=\sum_{j<i}H_{ij}u_j(z)$. Then $w_0=0$ and $w_1=1/b_1$.
For $i\ge1$, use $zu_j=b_ju_{j-1}+b_{j+1}u_{j+1}$, with the $j=0$ term
understood through $b_0=0$. The coefficient of $u_j$ in
$zw_i-b_iw_{i-1}-b_{i+1}w_{i+1}$ is
\[
 b_{j+1}H_{i,j+1}+b_jH_{i,j-1}
   -b_iH_{i-1,j}-b_{i+1}H_{i+1,j}.
\]
Entries with column $-1$ are read as zero. This coefficient
vanishes by triangularity when $j>i$, and by parity when $i-j$
is odd. For $j=i$ it is $1-1=0$.
If $i-j=2h\ge2$, put $s=(i+j)/2$. Substitution of
\eqref{eq:Hentries} reduces it to
\[
 (-1)^h\left(\frac{b_i+b_j}{b_s}
                 -\frac{b_{i+1}+b_{j+1}}{b_{s+1}}\right)=0,
\]
because both fractions equal $2\cos(\pi h\alpha)$.
This also covers $j=0$, since $b_0=0$; the denominators have
positive indices. Thus $w$ satisfies \eqref{eq:vrec} with
the initial values of $v$, and $w=v$. Triangularity gives
\begin{equation}\label{eq:vHp}
 v^{[n]}(z,\alpha)=H^{[n]}(\alpha)u^{[n]}(z,\alpha).
\end{equation}

{ (3) Multiply \eqref{eq:prec} by $v_j$ and \eqref{eq:vrec}
by $u_j$, and subtract. The result is $b_{j+1}W(u,v)(j+1)=b_jW(u,v)(j)$.
At $j=1$ the initial conditions give $b_1W(u,v)(1)=-1$.}
\end{proof}

\begin{remark}
Later, we will use certain Hankel structure hidden in $H(\alpha)$, so here we
present one concrete matrix $H^{[6]}(\alpha)$:
\[
 H^{[6]}(\alpha)=
 \begin{pmatrix}
 0&0&0&0&0&0\\
 b_1^{-1}&0&0&0&0&0\\
 0&b_2^{-1}&0&0&0&0\\
 -b_2^{-1}&0&b_3^{-1}&0&0&0\\
 0&-b_3^{-1}&0&b_4^{-1}&0&0\\
 b_3^{-1}&0&-b_4^{-1}&0&b_5^{-1}&0
 \end{pmatrix}.
\]
\end{remark}
\begin{lemma}\label{lem:globalH}
There exists a subset $\mathbb I_1\subset\mathbb I$ with
$|\mathbb I_1|=1$ such that, for every $\alpha\in\mathbb I_1$
and every positive integer $n$,
\[
 \norm{H^{[n]}(\alpha)}\lesssim_\alpha n^2.
\]
\end{lemma}
\begin{proof}
Write $\distT{x}=\min_{\ell\in\Z}|x-\ell|$ and define
\[
 \mathbb I_1=
 \bigcup_{m=1}^{\infty}\bigcap_{r=1}^{\infty}
 \{\alpha\in\mathbb I:\distT{r\alpha}\ge m^{-1}r^{-3/2}\}.
\]
For $\epsilon>0$,
$\Leb\{\alpha\in(0,1):\distT{r\alpha}<\epsilon\}\le2\epsilon$.
For each $m$, the complement of the corresponding intersection
therefore has measure at most
$2m^{-1}\sum_{r\ge1}r^{-3/2}\le6/m$. Taking $m\to\infty$ proves
$|\mathbb I_1|=1$. The sets used here are Borel.

Fix $\alpha\in\mathbb I_1$ and choose $a=a(\alpha)>0$ such that
$\distT{r\alpha}\ge ar^{-3/2}$ for all $r\ge1$.
The points $0,\alpha,\ldots,n\alpha$ are
$\epsilon_n$-separated modulo $1$, where $\epsilon_n=an^{-3/2}$.
On each side of $0$, the $j$th nearest point has distance at least
$j\epsilon_n$. Since $|b_r|\ge4\distT{r\alpha}$,
\[
 \sum_{r=1}^n|b_r|^{-1}
 \le\frac{1+\log n}{2a}n^{3/2}
 \le\frac{3}{2a}n^2.
\]
The last inequality follows from $1+\log n\le3\sqrt n$.
Since the denominator indices in a row or column of $H^{[n]}$ are
distinct elements of $\{1,\ldots,n-1\}$, all absolute row and column
sums of $H^{[n]}$ are bounded by this sum. For a finite real matrix $A$,
the standard row-column estimate \cite[Corollary~2.3.2, p.~73]{GVL2013} is
\begin{equation}\label{eq:row-column-bound}
 \begin{gathered}
 \norm A\le\sqrt{\norm A_1\norm A_\infty},\\
 \norm A_1=\max_j\sum_i|A_{ij}|,\qquad
 \norm A_\infty=\max_i\sum_j|A_{ij}|.
 \end{gathered}
\end{equation}
Thus {$\|H^{[n]}\|_1\le3n^2/(2a)$ and $\|H^{[n]}\|_\infty\le3n^2/(2a)$}.
Applying \eqref{eq:row-column-bound} to $H^{[n]}$ gives
\[
 \norm{H^{[n]}}\le\frac3{2a}n^2.
\]
\end{proof}

\begin{lemma}\label{lem:power-law}
There exist an absolute constant $\rho>0$ and a set
$\mathbb I_2\subset\mathbb I$ with $|\mathbb I_2|=1$ such that,
for every $\alpha\in\mathbb I_2$, $E\in\R$ and positive integer $n$,
\begin{equation}\label{eq:power-law}
 \norm{u^{[n]}(E,\alpha)}\gtrsim_\alpha n^\rho.
\end{equation}
\end{lemma}

We will prove Lemma~\ref{lem:power-law} in Section~\ref{sec:growth}.

\subsection{\texorpdfstring{$M$}{M}-function and proof of Lemma \ref{lem:holder}}

Given {$z\in\C\setminus\Sigma_\alpha$}, define 
\begin{equation}\label{eq:M-function}
  w(z,\alpha):=(J_\alpha^+-z)^{-1}\delta_0\ \ \text{ and }\ \  m(z,\alpha):=\langle \delta_0,(J_\alpha^+-z)^{-1}\delta_0\rangle=w_0(z,\alpha).
\end{equation}
By \eqref{eq:mu}, we have 
\begin{equation}\label{eq:m-integral}
  m(z,\alpha)=\int_{\Sigma_\alpha}\frac{1}{E'-z}d\mu_\alpha(E').
\end{equation}
{The difference $w-v$ solves the homogeneous
equation and has initial coordinate $w_0$.} By \eqref{eq:def-u-v}, we have 
\begin{equation}\label{eq:w-u-v}
  w(z,\alpha)=v(z,\alpha)+w_0(z,\alpha)u(z,\alpha)=v(z,\alpha)+m(z,\alpha)u(z,\alpha).
\end{equation}

\begin{proof}[Proof of Lemma \ref{lem:holder}]
  Fix $\alpha\in \I_1\cap \I_2$. In the following proof, we {suppress} $\alpha$ from the notation. Given  $z=E+r i$ with $E\in \R$ and $r>0$, then by \eqref{eq:m-integral},
\begin{equation}\label{eq:m-complex}
  m(z)=\int_{\Sigma_\alpha}
  \frac{E'-E+ri}{(E'-E)^2+r^2}d\mu_\alpha(E').
\end{equation} 
Since $\mu_\alpha$ is supported on $\Sigma_\alpha,$ we have 
\begin{align}\label{eq:mass-poisson}
\nonumber  r{\rm Im}\ m(z)&=\int_{\Sigma_\alpha}
  \frac{r^2}{(E'-E)^2+r^2}d\mu_\alpha(E')\ge \int_{\Sigma_\alpha\cap B(E,r)}
  \frac{r^2}{(E'-E)^2+r^2}d\mu_\alpha(E')\\
  &\ge\frac{\mu_\alpha(B(E,r))}{2}.
\end{align}

By \eqref{eq:mu} and \eqref{eq:m-complex}, we have 
\begin{align*}
  \|w(z)\|^2 & =\langle w(z),w(z)\rangle=\langle (J_\alpha^+-z)^{-1}\delta_0,(J_\alpha^+-z)^{-1}\delta_0\rangle\\
  &=\langle\delta_0, (J_\alpha^+-\bar z)^{-1}(J_\alpha^+-z)^{-1}\delta_0\rangle=\int_{\Sigma_\alpha}\frac{1}{(E'-\bar z)(E'-z)}d\mu_\alpha(E')\\
  &=\int_{\Sigma_\alpha}
  \frac{1}{(E'-E)^2+r^2}d\mu_\alpha(E')=\frac{{\rm Im }\ m(z)}{r}\le \frac{|m(z)|}{r}.
\end{align*}
 By \eqref{eq:w-u-v}, Lemma \ref{lem:solutions}(2) and Lemma \ref{lem:globalH}, there exists a constant $C_1=C_1(\alpha)>0$ such that 
\begin{align*}
  |m(z)|\|u^{[n]}(z)\| & =\|w^{[n]}(z)-v^{[n]}(z)\|=\|w^{[n]}(z)-H^{[n]}u^{[n]}(z)\|\\
  &\le\|w(z)\|+C_1n^2\|u^{[n]}(z)\|\le \sqrt{\frac{|m(z)|}{r}}+C_1n^2\|u^{[n]}(z)\|.
\end{align*}
  Hence we have 
  \begin{align*}
    |m(z)| & \le C_1n^2+\sqrt{\frac{|m(z)|}{r\|u^{[n]}(z)\|^2}}\le C_1n^2+\frac{|m(z)|}{2}+\frac{1}{2r\|u^{[n]}(z)\|^2}.
  \end{align*}
  By Lemma \ref{lem:solutions}(1), we have $\|u^{[n]}(z)\|\ge \|u^{[n]}(E)\|$. So by Lemma \ref{lem:power-law}, there exists a constant $C_2=C_2(\alpha)>0$ such that 
  \begin{align}\label{eq:rmbound}
    r|m(z)| & \le 2C_1rn^2+\|u^{[n]}(E)\|^{-2}\le 2C_1rn^2+C_2^{-2}n^{-2\rho}.
  \end{align}
Take $n=\lceil r^{-1/(2(\rho+1))}\rceil$ for $0<r<1$. Then
$r^{-1/(2(\rho+1))}\le n\le2r^{-1/(2(\rho+1))}$. Write $d:=\rho/(1+\rho)$.
Thus $rn^2\le4r^d$ and $n^{-2\rho}\le r^d$. By
\eqref{eq:mass-poisson} and \eqref{eq:rmbound},
\[
 \mu_\alpha(B(E,r))\le2r|m(z)|\le 2(8C_1+C_2^{-2})r^d .
\]
So the result follows.
\end{proof}
\section{Proof of Lemma~\ref{lem:power-law}}\label{sec:growth}

\subsection{A recursive lower bound for the norm of \texorpdfstring{$u^{[n]}(E,\alpha)$}{the canonical solution}}
\label{subsec:recursive}
Given $\alpha\in\I$, assume it has continued fraction expansion
$\alpha=[a_1(\alpha),a_2(\alpha),\ldots]$. Write
\begin{equation}\label{eq:A0}
 A_0:=2048.
\end{equation}
Define a set of indices as
\begin{equation}\label{eq:index-alpha}
 \mathcal I_\alpha:=\{n\in\N:\ {n\ge20,\ } a_{n+1}(\alpha)\ge A_0\}.
\end{equation}
Assume $p_n(\alpha)/q_n(\alpha)$ is the $n$th convergent of $\alpha$.
We have the following key recursive estimate.

\begin{lemma}\label{lem:recursive}
There is an absolute constant $D>0$ such that for any $\alpha\in\I$,
$E\in\R$ and $N\in\N$,
\[
 \sum_{\substack{n\in\mathcal I_\alpha\\q_n\le N/2}}
 \|u^{[{\lfloor q_n/8\rfloor}]}(E,\alpha)\|^2
 \le D\|u^{[N]}(E,\alpha)\|^2.
\]
\end{lemma}

We will prove Lemma~\ref{lem:recursive} in Section~\ref{sec:accumulation}.

\subsection{Power law estimate holds for good frequencies}\label{subsec:goodfrequencies}
In this part, we show that if $\alpha$ is good in the sense we will define
now, then, combined with Lemma~\ref{lem:recursive}, this gives the power
law estimate \eqref{eq:power-law}.

\begin{definition}\label{def:goodfrequencies}
Fix $\eta\in(0,1)$ and $\kappa\ge1$. Define two sets of frequencies as
\begin{align}
 \mathbb F_\eta
 &:=\{\alpha\in\I:\exists N(\alpha)\in\N\text{ s.t. }
 \#(\mathcal I_\alpha\cap[1,n])\ge\eta n,\ \forall n\ge N(\alpha)\},
 \label{eq:Feta}\\
 \mathbb G_\kappa
 &:=\{\alpha\in\I:\exists N(\alpha)\in\N\text{ s.t. }
 q_n(\alpha)\le e^{\kappa n},\ \forall n\ge N(\alpha)\}.
 \label{eq:Gkappa}
\end{align}
\end{definition}

At first we estimate the lower bounds of $\|u^{[n]}(E,\alpha)\|$ along a
subsequence for $\alpha\in\mathbb F_\eta$.

\begin{lemma}\label{lem:subsequence}
Assume $\alpha\in\mathbb F_\eta$ and $E\in\R$. Then
\[
 \|u^{[q_n]}(E,\alpha)\|\gtrsim_{\alpha,\eta}\theta(\eta)^n,\qquad
 \text{where }\quad
 \theta(\eta):=\bigl((D+1)/D\bigr)^{\eta/20}.
\]
\end{lemma}

\begin{proof}
Write $\widehat{\mathcal I}_\alpha:=\mathcal I_\alpha-8$.
For any $n\in\N$, define
\begin{equation}\label{eq:XY}
 X_n:=\|u^{[q_n]}(E,\alpha)\|^2
 \qquad\text{and}\qquad
 Y_n:=\sum_{j\in\widehat{\mathcal I}_\alpha\cap[1,n]}X_j.
\end{equation}
{Set $Y_0:=0$.}
If $\widehat{\mathcal I}_\alpha\cap[1,n]=\varnothing$, we set $Y_n:=0$.
Thus $Y_n$ is nondecreasing and
\begin{equation}\label{eq:XYdifference}
 X_n=Y_n-Y_{n-1}\qquad\text{if } n\in\widehat{\mathcal I}_\alpha.
\end{equation}
\noindent\textbf{Claim.} The following estimates hold:
\[
 \begin{cases}
  X_n\ge D^{-1}Y_{n-10},& n>10,\\
  Y_n\ge \dfrac{1+D}{D}\,Y_m,& n\in\widehat{\mathcal I}_\alpha,\ m\ge1,\ n-m\ge10.
 \end{cases}
\]
\noindent$\lhd$ It is well known that $q_{n+2}\ge2q_n$, thus
$16q_{m-8}\le q_m$ for $m\ge8$.
{For $m\in\mathcal I_\alpha$, the denominator recurrence gives
$q_m\ge m\ge20$, so $\lfloor q_m/8\rfloor\ge q_m/8-1\ge q_m/16
\ge q_{m-8}$. Also, $m\le n-2$ implies $q_m\le q_n/2$.}
By Lemma~\ref{lem:recursive},
\begin{align*}
 X_n&=\|u^{[q_n]}(E,\alpha)\|^2
 \ge\frac1D\sum_{\substack{m\in\mathcal I_\alpha\\q_m\le q_n/2}}
       \|u^{[\lfloor q_m/8\rfloor]}(E,\alpha)\|^2\\
 &\ge\frac1D\sum_{\substack{m\in\mathcal I_\alpha\\m\le n-2}}
       \|u^{[q_{m-8}]}(E,\alpha)\|^2
 =\frac1D\sum_{j\in\widehat{\mathcal I}_\alpha\cap[1,n-10]}X_j
 =\frac{Y_{n-10}}D .
\end{align*}
Now assume $n,m\in\N$ with $n\in\widehat{\mathcal I}_\alpha$ and
$n-m\ge10$. Since $Y_n$ is nondecreasing, by \eqref{eq:XYdifference},
\[
 Y_n=X_n+Y_{n-1}\ge D^{-1}Y_{n-10}+Y_m\ge(1+1/D)Y_m.
\]
So the Claim holds. \hfill$\rhd$

List the elements in $\widehat{\mathcal I}_\alpha\cap[1,n]$ as
\[
 m_0<m_1<\cdots<m_{k-1}<m_k .
\]
Since every $m\in\mathcal I_\alpha$ satisfies $m\ge20$, the map
$m\mapsto m-8$ sends $\mathcal I_\alpha\cap[1,n]$ into
$\widehat{\mathcal I}_\alpha\cap[1,n]$. Thus \eqref{eq:Feta} gives
$k+1\ge\eta n$ for all sufficiently large $n$.
Write $l=\lfloor k/10\rfloor$ and $\lambda:=(D+1)/D$.
{Also $Y_{m_0}=X_{m_0}\ge1$ since $u_0=1$,
and $m_{j+10}-m_j\ge10$.}
By the Claim, for $n$ large enough that $l\ge1$,
\begin{align*}
 X_n &\ge\frac{Y_{n-10}}D
 \ge\frac{Y_{m_k-10}}D\ge\frac{Y_{m_{10l}-10}}D
 \ge\frac{Y_{m_{10(l-1)}}}D\\
 &\ge\frac{\lambda^{l-1}}D
 \ge\frac{\lambda^{k/10-2}}D
 \ge\frac{\lambda^{\eta n/10}}{D\lambda^3}.
\end{align*}
{Taking square roots gives the exponent $\eta/20$ in the
bound for $\|u^{[q_n]}\|$. The finitely many smaller indices are covered
by $\|u^{[q_n]}\|\ge1$, uniformly in $E$.}
So the result follows.
\end{proof}

Next we obtain the power law estimate \eqref{eq:power-law} for
$\alpha\in\mathbb F_\eta\cap\mathbb G_\kappa$.

\begin{lemma}\label{lem:power-general}
Assume $\alpha\in\mathbb F_\eta\cap\mathbb G_\kappa$ and $E\in\R$.
Then
\[
 \|u^{[n]}(E,\alpha)\|\gtrsim_{\alpha,\eta,\kappa}
 n^{\rho(\eta,\kappa)},\qquad
 \text{where }\quad \rho(\eta,\kappa)=\kappa^{-1}\log\theta(\eta).
\]
\end{lemma}

\begin{proof}
Fix $\alpha\in\mathbb F_\eta\cap\mathbb G_\kappa$. Assume
$q_k\le n<q_{k+1}$. If $n$ is large enough, then $k$ is also large
enough. So by \eqref{eq:Gkappa}, we have $q_{k+1}\le e^{\kappa(k+1)}$.
By Lemma~\ref{lem:subsequence},
\begin{align*}
 \|u^{[n]}(E,\alpha)\|&\ge\|u^{[q_k]}(E,\alpha)\|
 \gtrsim_{\alpha,\eta}\theta(\eta)^k\\
 &\ge e^{(\kappa^{-1}\log q_{k+1}-1)\log\theta(\eta)}
 \ge\theta(\eta)^{-1}n^{\log\theta(\eta)/\kappa}.
\end{align*}
{The threshold for $n$ is independent of $E$. Once again
$\|u^{[n]}\|\ge1$ for the remaining lengths.}
\end{proof}

\subsection{Proof of Lemma~\ref{lem:power-law}}\label{subsec:fullmeasure-growth}
Now we apply the ergodicity of the Gauss measure to finish the proof of
Lemma~\ref{lem:power-law}.

\begin{proof}
It is classical that (see for example \cite[Corollary~3.8]{EW2011})
for a.e.\ $\alpha\in\I$,
\begin{equation}\label{eq:Levy}
 \frac{{\log} q_n(\alpha)}{n}\longrightarrow
 \frac{\pi^2}{12\log2}=1.18\ldots .
\end{equation}
Thus we have $|\mathbb G_2|=1$.

Define the Gauss map $G:\I\to\I$ by $G(\alpha):=\{1/\alpha\}$.
Define the Gauss measure on $[0,1]$ by
\[
 m_G(A):=\frac1{\log2}\int_A\frac{dx}{1+x}.
\]
Then $m_G$ is equivalent to Lebesgue measure on $[0,1]$.
It is well known (see for example \cite[Chapter~3]{EW2011}) that if
$\alpha=[a_1(\alpha),a_2(\alpha),\ldots]$, then
$G(\alpha)=[a_2(\alpha),a_3(\alpha),\ldots]$, and $m_G$ is ergodic
with respect to $G$. It is seen that $a_1(\alpha)\ge A_0$ if and only
if $\alpha\in(0,1/A_0)$, so by \eqref{eq:index-alpha}, for $n\ge20$,
\[
 \#(\mathcal I_\alpha\cap[1,n])
 =\sum_{i={20}}^n\chi_{(0,1/A_0)}(G^i(\alpha)).
\]
{Omitting the first nineteen terms does not change the limit
of the averages.} Then by ergodic theory, for $m_G$-a.e.\ $\alpha\in\I$,
\begin{align*}
 \lim_{n\to\infty}\frac1n\#(\mathcal I_\alpha\cap[1,n])
 &=\lim_{n\to\infty}\frac1n\sum_{i=1}^n\chi_{(0,1/A_0)}(G^i(\alpha))\\
 &=m_G((0,1/A_0))=\frac{\log(1+1/A_0)}{\log2}.
\end{align*}
Take
\begin{equation}\label{eq:eta}
 \eta:=\log(1+1/A_0).
\end{equation}
Then $\eta$ is an absolute constant. Since $\log2<1$ and $m_G$ is
equivalent to Lebesgue measure on $[0,1]$, we conclude that
$|\mathbb F_\eta|=1$. Define
\[
 \I_2:=\mathbb F_\eta\cap\mathbb G_2.
\]
Then $|\I_2|=1$ and by Lemma~\ref{lem:power-general}, for any
$\alpha\in\I_2$ we have
\[
 \|u^{[n]}(E,\alpha)\|\gtrsim_\alpha n^{\rho(\eta)},
\]
where
\[
 \rho(\eta)=\frac12\log\theta(\eta),\quad
 \theta(\eta)=\bigl((D+1)/D\bigr)^{\eta/20},\quad
 \eta=\log(1+1/A_0).
\]
Thus $\rho(\eta)$ is an absolute constant.
\end{proof}

\section{Proof of Lemma~\ref{lem:recursive}}\label{sec:accumulation}
Throughout this section, we fix $\alpha\in\I$ and $E\in\R$.
To simplify the notation, we write
\[
 u:=u(E,\alpha),\quad v:=v(E,\alpha),\quad H:=H(\alpha),\quad
 \mathcal I:=\mathcal I_\alpha,\quad q_k:=q_k(\alpha).
\]
The proof is based on two estimates: one is an upper bound for
$\|\Pi_{n_k}u\|$, {and the other is an upper bound for a
truncation operator.}

\subsection{Proof of Lemma~\ref{lem:recursive} by assuming two estimates}\label{sec:summation}
We fix the following notations:
\begin{equation}\label{eq:window}
 \begin{gathered}
 n_k:=\lfloor q_k/8\rfloor,\quad m_k:=q_k+n_k,\\
 I_k:=\{q_k-n_k,\ldots,m_k-1\},\qquad L_k:=\{n_k,\ldots,m_k-1\}.
 \end{gathered}
\end{equation}

\begin{lemma}\label{prop:area}
For any $k\in\mathcal I$, one has
\[
 \frac45\|\Pi_{n_k}u\|^2\le|b_{q_k}|\|\Pi_{I_k}u\|\|\Pi_{I_k}v\|.
\]
\end{lemma}

For any $k\in\N$, define the truncation operator
\begin{equation}\label{eq:Gamma-def}
 \Gamma_k:=b_{q_k}\Pi_{I_k}H^{[m_k]}\Pi_{m_k}.
\end{equation}
{Here $\Pi_{m_k}$ selects the first $m_k$ input coordinates and
$\Pi_{I_k}$ selects the output coordinates with indices in $I_k$.
When a restriction $\Pi_S\phi$ is used as the input of another
operator, we extend it by zero outside $S$, retaining the original
coordinate labels. In this sense $\Pi_{m_k}=\Pi_{n_k}+\Pi_{L_k}$.}

\begin{lemma}\label{lem:returnmatrix}
There exists an absolute constant $\varrho>0$ such that for any
$k\in\mathcal I$,
\[
 \|\Gamma_k\|\le\varrho.
\]
\end{lemma}

The proofs of these two lemmas are highly nontrivial, and will be given
in Section~\ref{sec:local} and Section~\ref{sec:returnnorm}, respectively.
Now let us assume them and finish the proof of Lemma~\ref{lem:recursive}.

We need one more observation. For any $N\in\N$, define
\begin{equation}\label{eq:returnindices}
 \mathcal I(N):=\{k\in\mathcal I:q_k\le N/2\}.
\end{equation}

\begin{lemma}\label{lem:summation}
For any $\phi\in\C^{\Z^+}$, we have
\begin{align*}
 \sum_{k\in\mathcal I(N)}\|\Pi_{I_k}\phi\|^2
 &\le\sum_{k\in\mathcal I(N)}\|\Pi_{L_k}\phi\|^2\le\|\Pi_N\phi\|^2,\\
 \biggl(\sum_{k\in\mathcal I(N)}\|\Gamma_k\phi\|^2\biggr)^{1/2}
 &\le\varrho\|\Pi_N\phi\|
   +\varrho\biggl(\sum_{k\in\mathcal I(N)}\|\Pi_{n_k}\phi\|^2\biggr)^{1/2}.
\end{align*}
\end{lemma}

\begin{proof}
At first, we have the following observation:

\noindent\textbf{Claim.}
\begin{enumerate}
 \item For any $k\in\mathcal I(N)$, $I_k\subset L_k\subset\{1,\ldots,N-1\}$.
 \item For any $k,l\in\mathcal I(N)$ with $k\ne l$,
       $L_k\cap L_l=\varnothing$. Hence $I_k\cap I_l=\varnothing$.
\end{enumerate}
\noindent$\lhd$ (1) By \eqref{eq:window} and \eqref{eq:returnindices},
{$q_k\ge k\ge20$, so that $n_k\ge2$ and
$n_k\ge q_k/8-1\ge q_k/16$. In particular,}
\[
 n_k\le\frac{q_k}8<\frac{7q_k}8\le q_k-n_k
 <m_k-1\le\frac{9q_k}8<2q_k\le N.
\]
So (1) follows.

(2) WLOG assume $k<l$. Since $k\in\mathcal I$,
\[
 q_l\ge q_{k+1}\ge a_{k+1}q_k\ge A_0q_k\ge32q_k.
\]
Hence
\[
 n_l\ge\frac{q_l}{16}\ge2q_k>m_k-1 .
\]
So $L_k\cap L_l=\varnothing$. Since $I_k\subset L_k$ and $I_l\subset L_l$,
$I_k\cap I_l=\varnothing$.\hfill$\rhd$

By the Claim,
\[
 I:=\bigsqcup_{k\in\mathcal I(N)}I_k
 \ \subset\ L:=\bigsqcup_{k\in\mathcal I(N)}L_k
 \ \subset\ \{0,\ldots,N-1\}.
\]
So we have
\[
 \sum_{k\in\mathcal I(N)}\|\Pi_{I_k}\phi\|^2
 =\|\Pi_I\phi\|^2\le\|\Pi_L\phi\|^2
 =\sum_{k\in\mathcal I(N)}\|\Pi_{L_k}\phi\|^2
 \le\|\Pi_N\phi\|^2.
\]
Thus the first inequality holds.

Since $\Pi_n$ is a projection, by \eqref{eq:Gamma-def} we have
\[
 \Gamma_k\phi=\Gamma_k\Pi_{m_k}\phi
 =\Gamma_k(\Pi_{n_k}+\Pi_{L_k})\phi
 =\Gamma_k\Pi_{n_k}\phi+\Gamma_k\Pi_{L_k}\phi.
\]
By Claim (2) and the definition of $\Gamma_k$, we can define three long vectors as
\[
 \psi:=(\Gamma_k\phi)_{k\in\mathcal I(N)},\quad
 \psi_1:=(\Gamma_k\Pi_{n_k}\phi)_{k\in\mathcal I(N)},\quad
 \psi_2:=(\Gamma_k\Pi_{L_k}\phi)_{k\in\mathcal I(N)}.
\]
{They belong to the same finite direct sum of Euclidean spaces,
with the sum of the squared block norms as its squared norm.}
Then $\psi=\psi_1+\psi_2$. So by Lemma~\ref{lem:returnmatrix} and Claim (1),
we have
\begin{align*}
 \biggl(\sum_{k\in\mathcal I(N)}\|\Gamma_k\phi\|^2\biggr)^{1/2}
 &=\|\psi\|\le\|\psi_1\|+\|\psi_2\|\\
 &=\biggl(\sum_{k\in\mathcal I(N)}\|\Gamma_k\Pi_{n_k}\phi\|^2\biggr)^{1/2}
 +\biggl(\sum_{k\in\mathcal I(N)}\|\Gamma_k\Pi_{L_k}\phi\|^2\biggr)^{1/2}\\
 &\le\varrho\biggl(\sum_{k\in\mathcal I(N)}\|\Pi_{n_k}\phi\|^2\biggr)^{1/2}
 +\varrho\biggl(\sum_{k\in\mathcal I(N)}\|\Pi_{L_k}\phi\|^2\biggr)^{1/2}\\
 &\le\varrho\biggl(\sum_{k\in\mathcal I(N)}\|\Pi_{n_k}\phi\|^2\biggr)^{1/2}
 +\varrho\|\Pi_N\phi\|.
\end{align*}
Thus the second inequality holds.
\end{proof}

Now we are ready for the proof.

\begin{proof}[Proof of Lemma~\ref{lem:recursive}]
Write
\[
 S_N:=\sum_{k\in\mathcal I(N)}\|\Pi_{n_k}u\|^2.
\]
By Lemma~\ref{lem:solutions}(2),
\[
 \Gamma_ku=b_{q_k}\Pi_{I_k}H^{[m_k]}\Pi_{m_k}u=b_{q_k}\Pi_{I_k}v.
\]
By Lemma~\ref{prop:area} and Lemma~\ref{lem:summation}, we have
\begin{align*}
 \frac45 S_N
 &\le\sum_{k\in\mathcal I(N)}|b_{q_k}|\|\Pi_{I_k}u\|\|\Pi_{I_k}v\|
 =\sum_{k\in\mathcal I(N)}\|\Pi_{I_k}u\|\|\Gamma_ku\|\\
 &\le\biggl(\sum_{k\in\mathcal I(N)}\|\Pi_{I_k}u\|^2\biggr)^{1/2}
      \biggl(\sum_{k\in\mathcal I(N)}\|\Gamma_ku\|^2\biggr)^{1/2}\\
 &\le\|\Pi_Nu\|\bigl(\varrho\|\Pi_Nu\|+\varrho\sqrt{S_N}\bigr)
 =\varrho\|\Pi_Nu\|^2+\varrho\|\Pi_Nu\|\sqrt{S_N}.
\end{align*}
Set $y:=\sqrt{S_N}/\|\Pi_Nu\|$ (here $\|\Pi_Nu\|\ge|u_0|=1$).
We obtain
\[
 \frac45y^2\le\varrho(1+y).
\]
Since $\varrho y\le\tfrac25 y^2+\tfrac58\varrho^2$, it follows
that $y^2\le\tfrac52\varrho+\tfrac{25}{16}\varrho^2$.
From this we get
\[
 \frac{S_N}{\|\Pi_Nu\|^2}=y^2\le5(1+\varrho)^2=:D .
\]
{Thus $S_N\le D\|\Pi_Nu\|^2$ for every $N\ge1$,
as required in Lemma~\ref{lem:recursive}.}
So the result follows.
\end{proof}

\section{Proof of Lemma~\ref{prop:area}}\label{sec:local}
Let us describe the rough idea of the proof. The key observation is the following.
If $\phi\in\C^{\Z}$ solves the equation $(J_\alpha-z)\phi=0$, then at any site
$n\in\N$, the one-sided vector $\hat\phi=(\phi_n,\phi_{n+1},\ldots)$ solves the
equation $(J^+_{\alpha,\theta}-z)\hat\phi=\hat c\delta_0$ for the phase $\theta=n\alpha$
and some number $\hat c$. Meanwhile, the signed reflection
$\tilde\phi:=((-1)^j\phi_{n-1-j})_{j\ge0}$ solves the equation
$(J^+_{\alpha,-\theta}-z)\tilde\phi=\tilde c\delta_0$
for the same $\theta=n\alpha$ and some number $\tilde c$.
The factor $(-1)^j$ compensates for the sign change
$b_j(-\alpha;\theta)=-b_j(\alpha;-\theta)$ in the reflected coefficients.
Now just like we did in Section~3.1,
we can define two special solutions $u(z,\alpha;\theta)$ and $v(z,\alpha;\theta)$
and express $\hat\phi$ as a linear combination of $u(z,\alpha;\theta)$ and $v(z,\alpha;\theta)$.
{All such expansions are taken on the finite segments where the required coefficients are nonzero.}
By the same reasoning as in Lemma~\ref{lem:solutions}(2), there exist two lower triangular
matrices $U(\alpha;\theta)$ and $V(\alpha;\theta)$ such that
$u(z,\alpha;\theta)=U(\alpha;\theta)u(z,\alpha)$ and
$v(z,\alpha;\theta)=V(\alpha;\theta)u(z,\alpha)$ on these segments.
Like this, we express $\hat\phi=W(\alpha;\theta)u(z,\alpha)$ for some matrix $W$.
Similarly, there exists some matrix $W'$ such that $\tilde\phi=W'u(z,\alpha)$
on the reflected segment.

Now back to our lemma. We need to bound $\|\Pi_{n_k}u\|$ from above by
$\|\Pi_{I_k}u\|$ and $\|\Pi_{I_k}v\|$. If we write
\[
 I_k^+:=\{q_k,\ldots,m_k-1\},\quad I_k^-:=\{q_k-n_k,\ldots,q_k-1\},
 \quad \hat u^+:=\Pi_{I_k^+}u,\quad \hat u^-:=\Pi_{I_k^-}u,
\]
then $\Pi_{I_k}u=\hat u^-\hat u^+$, where juxtaposition denotes concatenation
in the original coordinate order.
For a vector $x=(x_0,\ldots,x_{n-1})$, define its signed reflection by
\[
 x^\sharp:=(x_{n-1},-x_{n-2},\ldots,(-1)^{n-1}x_0).
\]
Apply the above observation to $\phi=u$.
Then there exist two matrices $W^+$ and $W^-$ such that
\[
 \hat u^+=W^+\Pi_{n_k}u,\qquad (\hat u^-)^\sharp=W^-\Pi_{n_k}u,
\]
This establishes an exact relation between $\Pi_{I_k}u$ and $\Pi_{n_k}u$.
A similar relation holds for $\Pi_{I_k}v$ and $\Pi_{n_k}u$.
Now by exploring these two relations, we are able to figure out the proof of
Lemma~\ref{prop:area}.

\subsection{A general estimate}
\subsubsection{Two special formal solutions related to fixed phase}
Fix a phase $\theta$ {with its chosen real representative}. Define
$u(z,\alpha;\theta),v(z,\alpha;\theta)$ as follows:
\begin{equation}\label{eq:six-canonical}
 \begin{cases}
 (J^+_{\alpha,\theta}-z)u(z,\alpha;\theta)=0, & u_0(z,\alpha;\theta)=1,\\
 (J^+_{\alpha,\theta}-z)v(z,\alpha;\theta)=\delta_0, & v_0(z,\alpha;\theta)=0 .
 \end{cases}
\end{equation}
{More precisely, their first $n$ coordinates are defined if
$b_s(\alpha;\theta)\ne0$ for $1\le s<n$. For $n\ge2$ the equations at $0$
give $u_1=z/b_1(\alpha;\theta)$ and $v_1=1/b_1(\alpha;\theta)$.
The interior equations through $n-2$ determine the remaining coordinates.
For $n=1$ only the prescribed zeroth coordinates are used.
{We call $n\in\N$ admissible for $\theta$ if
$b_s(\alpha;\theta)b_s(\alpha;-\theta)\ne0$ for $1\le s<n$.
This condition will be used whenever both halves of a window of length $2n$ occur.}}
Thus $u(z,\alpha;0)=u(z,\alpha)$ and $v(z,\alpha;0)=v(z,\alpha)$, where
$u(z,\alpha)$ and $v(z,\alpha)$ are defined in \eqref{eq:def-u-v}.

For a sequence $y=(\ldots,y_{-2},y_{-1})$ on the negative half-line,
the same signed reflection is written as
\begin{equation}\label{eq:six-reflection}
 (y^\sharp)_j:=(-1)^jy_{-1-j}\qquad(j\ge0).
\end{equation}
On a finite segment $\{-n,\ldots,-1\}$ this formula is used for $0\le j<n$.
For any $\phi\in\C^\Z$, define $\phi^+:=\Pi_{\Z^+}\phi$ and
$\phi^-:=\Pi_{-\N}\phi$.
For a vector on $J_n:=\{-n,\ldots,n-1\}$ we use the same two
restrictions, with signed reflection on its negative indices.

\begin{lemma}\label{lem:six-connection}
{For every $n\in\N$, there exist two families of lower triangular
$n\times n$ matrices $U(\alpha;\theta),V(\alpha;\theta)$, defined for all phases
$\theta$ for which $n$ is admissible, such that the following hold.}
\begin{enumerate}
\item $U(\alpha;0)=I$ and $V(\alpha;0)=H(\alpha)$ on the segment. Moreover,
with all vectors restricted to their first $n$ coordinates,
\[
 \begin{aligned}
 u(z,\alpha;\theta)&=U(\alpha;\theta)u(z,\alpha),&
 v(z,\alpha;\theta)&=V(\alpha;\theta)u(z,\alpha),\\
 {u(z,\alpha;-\theta)}&={U(\alpha;-\theta)}u(z,\alpha),&
 {v(z,\alpha;-\theta)}&={V(\alpha;-\theta)}u(z,\alpha).
 \end{aligned}
\]
\item If $\phi$ solves $(J_{\alpha,\theta}-z)\phi=0$
{at the interior sites $-n+1,\ldots,n-2$ of $J_n$}, then, on
the same segment and with $b_0=b_0(\alpha;\theta)$,
\[
 \begin{aligned}
 \phi^+&=\phi_0u(z,\alpha;\theta)-\phi_{-1}b_0v(z,\alpha;\theta)\\
        &=[\phi_0U(\alpha;\theta)-\phi_{-1}b_0 V(\alpha;\theta)]u(z,\alpha),\\
 (\phi^-)^\sharp
       &=\phi_{-1}{u(z,\alpha;-\theta)}
            -\phi_0b_0{v(z,\alpha;-\theta)}\\
       &=[-\phi_0b_0{V(\alpha;-\theta)}
                 +\phi_{-1}{U(\alpha;-\theta)}]u(z,\alpha).
 \end{aligned}
\]
\item For any such $n\in\N$, we have
\[
 \binom{\Pi_n\phi^+}{\Pi_n(\phi^-)^\sharp}
   =\Omega_n(\alpha;\theta)
      \binom{\phi_0\Pi_nu(z,\alpha)}{\phi_{-1}\Pi_nu(z,\alpha)},
\]
where $\Omega_n(\alpha;\theta)$ is defined as
\begin{equation}\label{eq:six-omega}
 \Omega_n(\alpha;\theta):=
 \begin{pmatrix}
 U^{[n]}(\alpha;\theta)&-b_0(\alpha;\theta)V^{[n]}(\alpha;\theta)\\
 -b_0(\alpha;\theta){V^{[n]}(\alpha;-\theta)}&
 {U^{[n]}(\alpha;-\theta)}
 \end{pmatrix}.
\end{equation}
\end{enumerate}
\end{lemma}
{
\begin{proof}
By the recurrence, $u_j(z,\alpha)$ is a polynomial of degree $j$ with
leading coefficient $(b_1(\alpha)\cdots b_j(\alpha))^{-1}$.
The first $n$ of these polynomials thus form a basis for the polynomials
of degree less than $n$. The shifted $u_j$ and $v_j$ have degrees $j$ and
$j-1$, respectively, with $v_0=0$.
{For each $\gamma\in\{\theta,-\theta\}$, expansion in this
basis uniquely defines $U(\alpha;\gamma)$ and $V(\alpha;\gamma)$ and gives (1).
The matrix $V(\alpha;\gamma)$ is strictly lower triangular.}
At phase zero the identities follow from uniqueness of the expansion and
Lemma~\ref{lem:solutions}(2).

For (2), the positive half satisfies the phase-$\theta$ interior recurrence.
At the left endpoint $b_1(\alpha;\theta)\phi_1=z\phi_0-b_0\phi_{-1}$.
On the negative half put $y_j=(-1)^j\phi_{-1-j}$. The identity
$b_{-j}(\alpha;\theta)=-b_j(\alpha;-\theta)$ yields the phase-$(-\theta)$
interior recurrence and
$b_1(\alpha;-\theta)y_1=z\phi_{-1}-b_0\phi_0$. The initial values and
uniqueness now give the stated expansions. For $n=1$ the expansions
just give $\phi_0$ and $\phi_{-1}$. Writing the two expansions as a
single matrix equation gives (3).
\end{proof}}

The lemma has the following consequence.
\begin{corollary}\label{cor:six-ineq}
Assume $\Omega_n(\alpha;\theta)$ is invertible {for some $n\in\N$ admissible for $\theta$}.
Define a matrix
\begin{equation}\label{eq:six-delta}
 \Delta_n(\alpha;\theta):=(\Omega_n(\alpha;\theta)^{-1})^T
   \begin{pmatrix}0&-I_n\\ I_n&0\end{pmatrix}\Omega_n(\alpha;\theta)^{-1}.
\end{equation}
Then for any two solutions $\phi,\psi$ of the eigen-equation
$(J_{\alpha,\theta}-E)x=0$ {on the interior of $J_n$}, one has
\begin{equation}\label{eq:six-w-bound}
 |W(\phi,\psi)(0)|\|\Pi_nu(E,\alpha)\|^2
 \le \|\Pi_{J_n}\phi\|\,\|\Delta_n(\alpha;\theta)\|\,\|\Pi_{J_n}\psi\|.
\end{equation}
\end{corollary}
\begin{proof}
Since $\Omega_n(\alpha;\theta)$ is invertible, by Lemma~\ref{lem:six-connection}(3),
writing $P=\Pi_nu(E,\alpha)$ and $Y_\phi=(\Pi_n\phi^+,\Pi_n(\phi^-)^\sharp)^T$
as a block column, we have
\[
 \binom{\phi_0P}{\phi_{-1}P}=\Omega_n(\alpha;\theta)^{-1}Y_\phi,\qquad
 \binom{\psi_0P}{\psi_{-1}P}=\Omega_n(\alpha;\theta)^{-1}Y_\psi .
\]
So for $E\in\R$, we have
\[
 -W(\phi,\psi)(0)\|P\|^2
 =\binom{\phi_0P}{\phi_{-1}P}^{T}
 \begin{pmatrix}0&-I_n\\I_n&0\end{pmatrix}\binom{\psi_0P}{\psi_{-1}P}
 =Y_\phi^T\Delta_n(\alpha;\theta)Y_\psi.
\]
Notice that $Y_\phi$ is obtained from $\Pi_{J_n}\phi$ by reordering
the coordinates and {multiplying them by signs}. So the two vectors
have the same norm. We then have \eqref{eq:six-w-bound}.
\end{proof}

\subsubsection{Criterion for the invertibility of \texorpdfstring{$\Omega_n(\alpha;\theta)$}{Omega}}
In this part, we present a sufficient condition such that $\Omega_n(\alpha;\theta)$
is invertible. For any $n\in\N$, define
\begin{equation}\label{eq:six-lambda}
 \Lambda_n(\alpha;\theta):=\Omega_n(\alpha;\theta)\mathbb J_n\Omega_n(\alpha;\theta)^T,
 \qquad \mathbb J_n:=\begin{pmatrix}0&I_n\\-I_n&0\end{pmatrix}.
\end{equation}
We note that if $\Omega_n(\alpha;\theta)$ is invertible, then
$\Lambda_n(\alpha;\theta)=\Delta_n(\alpha;\theta)^{-1}$ by \eqref{eq:six-delta}.
So we only need to study when $\Lambda_n(\alpha;\theta)$ is invertible.
At $\theta=0$, we have $\Omega_n(\alpha;0)=I_{2n}$ and hence
$\Lambda_n(\alpha;0)=\mathbb J_n$. The following lemma estimates the
deviation from this identity for small $\theta$.

\begin{lemma}\label{lem:six-criterion}
Fix $\delta>0$ small. Assume $n\in\N$ is such that
\begin{equation}\label{eq:six-separation}
 \|r\alpha\|_{\T}\ge\delta \qquad (r\in\Z,\ 0<|r|<n).
\end{equation}
If $|\theta|<\delta/200$, then
\begin{equation}\label{eq:six-lambda-bound}
 \|\Lambda_n(\alpha;\theta)-\mathbb J_n\|\le\frac15.
\end{equation}
Consequently, $\Omega_n(\alpha;\theta)$ is invertible and
\begin{equation}\label{eq:six-inverse-bound}
 \|\Delta_n(\alpha;\theta)\|\le\frac54 .
\end{equation}
\end{lemma}
The proof of Lemma~\ref{lem:six-criterion} is a bit technical and will be given in
{Section~\ref{subsec:six-criterion-proof}}.

\subsection{Proof of Lemma~\ref{prop:area}}
Fix $\alpha\in\I$ and $E\in\R$. Let $u(E,\alpha)$ and $v(E,\alpha)$ be defined
by \eqref{eq:def-u-v}. Both solve the homogeneous interior recurrence at
every site $j\ge1$. We use this recurrence only within the return window.
Now fix $k\in\mathcal I_\alpha$ and put $\sigma_k:=(-1)^{p_k}$.
Define $\phi,\psi$ on $J_{n_k}$ as
\begin{equation}\label{eq:six-shifts}
 \phi_j:=\sigma_k^j u_{q_k+j}(E,\alpha),\qquad
 \psi_j:=\sigma_k^j v_{q_k+j}(E,\alpha),\quad -n_k\le j<n_k.
\end{equation}
Define a phase
\begin{equation}\label{eq:six-phase}
 \theta_k=q_k\alpha-p_k .
\end{equation}
It is well-known that
\begin{equation}\label{eq:convergent-bounds}
 \frac1{2q_kq_{k+1}}\le\left|\alpha-\frac{p_k}{q_k}\right|
 \le\frac1{q_kq_{k+1}} .
\end{equation}
Since $k\in\mathcal I_\alpha$, we have $a_{k+1}\ge A_0$. So
\begin{equation}\label{eq:six-phase-bound}
 |\theta_k|\le \frac1{q_{k+1}}\le\frac1{a_{k+1}q_k}\le\frac1{A_0q_k}.
\end{equation}

\begin{lemma}\label{lem:six-shifts}
$\phi$ and $\psi$ are solutions of $(J_{\alpha,\theta_k}-E)x=0$
on the interior of $J_{n_k}$, with
\begin{equation}\label{eq:six-w}
 |W(\psi,\phi)(0)|=\frac1{|b_{q_k}(\alpha)|}.
\end{equation}
\end{lemma}
{
\begin{proof}
By \eqref{eq:bcoeff}, $b_{q_k+j}(\alpha)=\sigma_k b_j(\alpha;\theta_k)$.
Inserting $w_{q_k+j}=\sigma_k^j x_j$ in the recurrence for $w=u$ or $v$
gives the asserted equation for $x$ throughout the interior, since
$q_k-n_k+1\ge1$. By the definitions of $\phi$ and $\psi$ and
Lemma~\ref{lem:solutions}(3), we have $W(\psi,\phi)(0)=\sigma_k W(v,u)(q_k)
=\sigma_k/b_{q_k}(\alpha)$.
\end{proof}}

Now we are ready for the proof of Lemma~\ref{prop:area}.
\begin{proof}[Proof of Lemma~\ref{prop:area}]
By \eqref{eq:six-shifts}, we have
\begin{equation}\label{eq:six-window}
 {\|\Pi_{J_{n_k}}\phi\|=\|\Pi_{I_k}u(E,\alpha)\|,
 \qquad \|\Pi_{J_{n_k}}\psi\|=\|\Pi_{I_k}v(E,\alpha)\|.}
\end{equation}
By the theory of continued fractions \cite[Chapters~I--II]{Khinchin},
\begin{equation}\label{eq:six-dk}
 \|r\alpha\|_\T\ge|q_{k-1}\alpha-p_{k-1}|=:\delta_k,\qquad
 r\in\Z,\quad 0<|r|<q_k .
\end{equation}
So \eqref{eq:six-separation} holds for $\delta=\delta_k$ and $n=n_k$.
By \eqref{eq:convergent-bounds}, $1/(2q_k)\le\delta_k\le1/q_k$.
So by \eqref{eq:six-phase-bound}, $|\theta_k|\le 2\delta_k/A_0<\delta_k/200$.
Then by Lemma~\ref{lem:six-criterion}, we have
\[
 \|\Delta_{n_k}(\alpha;\theta_k)\|\le \frac54 .
\]
Now by \eqref{eq:six-w-bound}, \eqref{eq:six-w} and \eqref{eq:six-window},
\[
 \|\Pi_{n_k}u(E,\alpha)\|^2
 \le \frac{\|\Pi_{J_{n_k}}\phi\|\,\|\Delta_{n_k}(\alpha;\theta_k)\|\,
             \|\Pi_{J_{n_k}}\psi\|}{|W(\phi,\psi)(0)|}
 \le \frac54 |b_{q_k}(\alpha)|\,\|\Pi_{I_k}u\|\,\|\Pi_{I_k}v\|.
\]
Thus the lemma holds.
\end{proof}

{
\subsection{Proof of Lemma~\ref{lem:six-criterion}}\label{subsec:six-criterion-proof}
\subsubsection{The structure of \texorpdfstring{$\Lambda_n(\alpha;\theta)$}{Lambda}}
We start with an explicit formula for $\Lambda_n(\alpha;\theta)$.
Fix $\alpha\in\I$. For any phase $\theta$ and any $n$ admissible for $\theta$, define
\begin{equation}\label{eq:six-CQ}
\begin{aligned}
 C_n(\alpha;\theta)&:=U^{[n]}(\alpha;\theta)V^{[n]}(\alpha;\theta)^T
                  -V^{[n]}(\alpha;\theta)U^{[n]}(\alpha;\theta)^T,\\
 Q_n(\alpha;\theta)&:=U^{[n]}(\alpha;\theta)U^{[n]}(\alpha;-\theta)^T
                  -b_0(\alpha;\theta)^2V^{[n]}(\alpha;\theta)V^{[n]}(\alpha;-\theta)^T.
\end{aligned}
\end{equation}
By \eqref{eq:six-omega} and \eqref{eq:six-lambda}, we have
\begin{equation}\label{eq:six-CQ-blocks}
 \Lambda_n(\alpha;\theta)=
 \begin{pmatrix}
 -b_0(\alpha;\theta)C_n(\alpha;\theta)&Q_n(\alpha;\theta)\\
 -Q_n(\alpha;\theta)^T&b_0(\alpha;\theta)C_n(\alpha;-\theta)
 \end{pmatrix}.
\end{equation}
We have the following explicit formulas for $C_n(\alpha;\theta)$ and
$Q_n(\alpha;\theta)$, where $\kappa_{ij}$ and $\sigma_{ij}$ are as in
\eqref{eq:matrix-indices}.
\begin{lemma}\label{lem:six-entries}
Fix $n\in\N$. Assume $b_s(\alpha;\theta)b_s(\alpha;-\theta)\ne0$ for
$1\le s<n$. Then for any $0\le i,j<n$,
\begin{equation}\label{eq:six-C-entries}
 (C_n(\alpha;\theta))_{ij}=
 \begin{cases}
 \dfrac{(-1)^{\sigma_{ij}+1}}{b_{\kappa_{ij}}(\alpha;\theta)},&i-j\text{ odd},\\[4pt]
 0,&i-j\text{ even},
 \end{cases}
\end{equation}
and
\begin{equation}\label{eq:six-Q-entries}
 (Q_n(\alpha;\theta))_{ij}=
 \begin{cases}
 1,&i=j,\\[2pt]
 \dfrac{(-1)^k b_0(\alpha;\theta)}{b_k(\alpha;\theta)},&i-j=2k\ne0,\\[4pt]
 0,&i-j\text{ odd}.
 \end{cases}
\end{equation}
\end{lemma}
The proof of this lemma is not hard, but quite tedious. So we will postpone it
to Section~\ref{sec:six-tech}. Now we use this lemma to study the invertibility
of $\Lambda_n(\alpha;\theta)$.

\subsubsection{The idea of the proof for Lemma~\ref{lem:six-criterion}}
Let us explain the main strategy of the proof. By \eqref{eq:six-CQ-blocks},
to estimate the norm of $\Lambda_n(\alpha;\theta)-\mathbb J_n$, we only
need to estimate the norms of $C_n(\alpha;\pm\theta)$ and $Q_n(\alpha;\theta)-I_n$.
To gain some intuition, let us take $n=7$ and write down the concrete form
of $C_7$ and $Q_7$. For this example, write $c_k=1/b_k(\alpha;\theta)$
and $q_k=1/b_k(\alpha;-\theta)$ for $0<|k|<7$.

At first we deal with $C_7$. By \eqref{eq:six-C-entries}, we have
\[
 C_7=\begin{pmatrix}
 0&c_1&0&-c_2&0&c_3&0\\
 -c_1&0&c_2&0&-c_3&0&c_4\\
 0&-c_2&0&c_3&0&-c_4&0\\
 c_2&0&-c_3&0&c_4&0&-c_5\\
 0&c_3&0&-c_4&0&c_5&0\\
 -c_3&0&c_4&0&-c_5&0&c_6\\
 0&-c_4&0&c_5&0&-c_6&0
 \end{pmatrix}.
\]
Let $P$ be the permutation matrix that puts the coordinates in the order
$(0,2,4,6,1,3,5)$. Then
\[
 P C_7 P^T=
 \begin{pmatrix}0&\widehat C_7^{\rm eo}\\
               -(\widehat C_7^{\rm eo})^T&0\end{pmatrix},
 \qquad
 \widehat C_7^{\rm eo}=
 \begin{pmatrix}
 c_1&-c_2&c_3\\
 -c_2&c_3&-c_4\\
 c_3&-c_4&c_5\\
 -c_4&c_5&-c_6
 \end{pmatrix}.
\]
Here $\rm eo$ refers to even rows and odd columns in the original order.
Let $D_1=\diag(1,-1,1,-1)$ and $D_2=\diag(1,-1,1)$. Removing the signs gives
\[
 C_7^{\rm eo}:=D_1\widehat C_7^{\rm eo}D_2=
 \begin{pmatrix}
 c_1&c_2&c_3\\
 c_2&c_3&c_4\\
 c_3&c_4&c_5\\
 c_4&c_5&c_6
 \end{pmatrix}.
\]
Thus $\|C_7\|=\|C_7^{\rm eo}\|$. Now observe that
\begin{equation}\label{eq:six-example-kernel}
 (C_7^{\rm eo})_{ij}=c_{i+j+1}
 =\frac{1}{2\sin\pi(x_i-y_j)},\qquad
 x_i=(i+1)\alpha+\theta,\quad y_j=-j\alpha .
\end{equation}
This enables us to use a result of Montgomery and Vaughan to estimate
the norm of $C_7^{\rm eo}$ (see Lemma~\ref{lem:six-cosecant} below).

Next we deal with $Q_7$. It is convenient to display the transpose,
which has the same norm. Since $b_{-k}(\alpha;\theta)=-b_k(\alpha;-\theta)$,
\eqref{eq:six-Q-entries} gives
\[
 (Q_7-I_7)^T=b_0(\alpha;\theta)
 \begin{pmatrix}
 0&0&q_{-1}&0&-q_{-2}&0&q_{-3}\\
 0&0&0&q_{-1}&0&-q_{-2}&0\\
 q_1&0&0&0&q_{-1}&0&-q_{-2}\\
 0&q_1&0&0&0&q_{-1}&0\\
 -q_2&0&q_1&0&0&0&q_{-1}\\
 0&-q_2&0&q_1&0&0&0\\
 q_3&0&-q_2&0&q_1&0&0
 \end{pmatrix}.
\]
The same permutation now places the even-even and odd-odd blocks on
the diagonal. With $S=\diag(D_1,D_2)$, we obtain
\[
 S P(Q_7-I_7)^T P^T S
 =-b_0(\alpha;\theta)\begin{pmatrix}Q_7^{\rm ee}&0\\0&(Q_7^{\rm ee})^{[3]}\end{pmatrix},
 \qquad
 Q_7^{\rm ee}=
 \begin{pmatrix}
 0&q_{-1}&q_{-2}&q_{-3}\\
 q_1&0&q_{-1}&q_{-2}\\
 q_2&q_1&0&q_{-1}\\
 q_3&q_2&q_1&0
 \end{pmatrix}.
\]
The odd-odd block is a restriction of the even-even block, so
$\|Q_7-I_7\|=|b_0(\alpha;\theta)|\,\|Q_7^{\rm ee}\|$.
For $i\ne j$ its entries are
$(Q_7^{\rm ee})_{ij}=[2\sin\pi((i-j)\alpha-\theta)]^{-1}$.
At $\theta=0$ this is the cosecant matrix on the nodes $j\alpha$ with
zero diagonal. Since $\theta$ is small, the idea is to treat
$Q_7^{\rm ee}$ as a perturbation of this matrix. Now we give the detailed proof.

\subsubsection{The structures of \texorpdfstring{$C_n(\alpha;\theta)$ and $Q_n(\alpha;\theta)$}{C and Q}}
Fix $n\in\N$ and assume $b_k(\alpha;\theta)b_k(\alpha;-\theta)\ne0$ for
$1\le k<n$. For $0<|k|<n$ define
\begin{equation}\label{eq:six-cq-sequences}
 \mathbf c_k:=\frac{1}{b_k(\alpha;\theta)},\qquad
 \mathbf q_k:=\frac{1}{b_k(\alpha;-\theta)},\qquad \mathbf q_0:=0.
\end{equation}
Write $n^*:=\lceil n/2\rceil$ and $n_*:=\lfloor n/2\rfloor$.
There are $n^*$ even labels and $n_*$ odd labels in $\{0,\ldots,n-1\}$.
Define $C_n^{\rm eo}\in\R^{n^*\times n_*}$ and
$Q_n^{\rm ee}\in\R^{n^*\times n^*}$ by
\begin{equation}\label{eq:six-reduced-entries}
 (C_n^{\rm eo})_{ij}:=\mathbf c_{i+j+1},\qquad
 (Q_n^{\rm ee})_{ij}:=\mathbf q_{i-j},
\end{equation}
with indices starting at $0$. We have the following block structure.

\begin{lemma}\label{lem:six-parity}
There exist unitary matrices $U_1,U_2,V_1,V_2$ such that
\[
 \begin{aligned}
 C_n(\alpha;\theta)&=U_1
 \begin{pmatrix}0&C_n^{\rm eo}\\(C_n^{\rm eo})^T&0\end{pmatrix}U_2,\\
 Q_n(\alpha;\theta)-I_n&=b_0(\alpha;\theta)V_1
 \begin{pmatrix}Q_n^{\rm ee}&0\\0&(Q_n^{\rm ee})^{[n_*]}\end{pmatrix}V_2 .
 \end{aligned}
\]
\end{lemma}
\begin{proof}
Put the even coordinates first, followed by the odd coordinates, each in
increasing order. In the even-row, odd-column block of $C_n$, the entry
at $(2i,2j+1)$ is $(-1)^{i-j}\mathbf c_{i+j+1}$ by \eqref{eq:six-C-entries}.
Multiplying row $i$ and column $j$ of this block by $(-1)^i$ and $(-1)^j$
respectively gives $C_n^{\rm eo}$. The opposite block is its negative
transpose by skew-symmetry. Changing the signs of the odd output rows
removes this last minus sign.

For $Q_n-I_n$, entries of opposite parities vanish. In its transpose,
the entry at $(2i,2j)$ or $(2i+1,2j+1)$, with $i\ne j$, is
$(-1)^{i-j+1}b_0(\alpha;\theta)\mathbf q_{i-j}$ by \eqref{eq:six-Q-entries}.
Both diagonal blocks have zero diagonal. The same alternating signs,
followed by an overall minus sign, give $b_0(\alpha;\theta)Q_n^{\rm ee}$
and its leading $n_*\times n_*$ block. Transposing back transposes
each block. Reversing its row and column orders restores the original
block, since its entries depend only on $i-j$. All operations used on
the left and right are unitary. For $n=1$ the odd block is empty and
both asserted identities reduce to zero.
\end{proof}

\subsubsection{Estimates of the norms of \texorpdfstring{$C_n(\alpha;\theta)$ and $Q_n(\alpha;\theta)-I_n$}{C and Q minus I}}
The main tool is the following consequence of the Montgomery--Vaughan
inequality.

\begin{lemma}[Cosecant kernels]\label{lem:six-cosecant}
Let $x_1,\ldots,x_m\in\mathbb R$ satisfy
$\|x_i-x_j\|_{\mathbb T}\ge\varepsilon>0$ for $i\ne j$. Define
\[
 \cZ_{ij}:=\frac{1}{2\sin\pi(x_i-x_j)}\quad(i\ne j),
 \qquad \cZ_{ii}:=0.
\]
Then
\[
 \|\cZ\|\le\frac{1}{2\varepsilon}.
\]
For any two subsets of the nodes, the submatrix obtained by selecting
the first subset as rows and the second as columns satisfies the same bound.
\end{lemma}
\begin{proof}
The case $m=1$ is immediate. For $m\ge2$, the cosecant form of the
Montgomery--Vaughan inequality \cite[Theorem~1, equation~(1.2)]{MV} gives
\[
 \left|\sum_{i\ne j}a_i\overline{a_j}\,
     \csc\bigl(\pi(x_i-x_j)\bigr)\right|
 \le\frac{1}{\varepsilon}\sum_i|a_i|^2
 \qquad(a\in\mathbb C^m).
\]
Since $\cZ$ is real and skew-symmetric, $\mathrm i\cZ$ is Hermitian.
The displayed inequality bounds the absolute quadratic form of
$2\mathrm i\cZ$ on unit vectors by $\varepsilon^{-1}$, and hence
$\|\cZ\|\le1/(2\varepsilon)$. Restricting rows and columns cannot
increase the norm, which gives the rectangular assertion.
Changing a node by an integer multiplies its row and column by signs,
so the choice of representatives modulo $1$ does not affect the bound.
\end{proof}

\Needspace{5\baselineskip}
Now we have the following estimates.
\begin{lemma}\label{lem:six-block-bounds}
Assume $n\in\N$ and $\delta>0$ satisfy \eqref{eq:six-separation}, and
$|\theta|\le\delta/200$. Then
\[
 |b_0(\alpha;\theta)|\,\|C_n(\alpha;\pm\theta)\|\le\frac{10|\theta|}{\delta},
 \qquad
 \|Q_n(\alpha;\theta)-I_n\|\le\frac{10|\theta|}{\delta}.
\]
\end{lemma}
\begin{proof}
For $n=1$ both quantities vanish. Suppose $n\ge2$. For $1\le s<n$ we have
$\|s\alpha\pm\theta\|_{\T}\ge\delta-|\theta|>0$, so both phases are
admissible. By Lemma~\ref{lem:six-parity},
\begin{equation}\label{eq:six-reduced-norms}
 \|C_n(\alpha;\theta)\|=\|C_n^{\rm eo}\|,\qquad
 \|Q_n(\alpha;\theta)-I_n\|=|b_0(\alpha;\theta)|\,\|Q_n^{\rm ee}\|.
\end{equation}
Indeed, the off-diagonal blocks of the first matrix are transposes of
each other, and the smaller diagonal block of the second is a compression
of the larger one.

As in \eqref{eq:six-example-kernel}, define $n$ nodes by
\begin{equation}\label{eq:six-kernel-nodes}
 t_i:=\begin{cases}
 (i+1)\alpha+\theta,&0\le i<n^*,\\
 -(i-n^*)\alpha,&n^*\le i<n.
 \end{cases}
\end{equation}
Then $(C_n^{\rm eo})_{ij}=[2\sin\pi(t_i-t_{n^*+j})]^{-1}$.
Within either group of nodes, \eqref{eq:six-separation} gives separation
at least $\delta$ modulo $1$. For a pair from different groups,
$1\le i+j+1\le n-1$, so
\[
 \|t_i-t_{n^*+j}\|_{\T}
 \ge\|(i+j+1)\alpha\|_{\T}-|\theta|\ge\delta-|\theta|>\delta/2 .
\]
Applying Lemma~\ref{lem:six-cosecant} to these row and column sets yields
\begin{equation}\label{eq:six-c-reduced-bound}
 \|C_n^{\rm eo}\|\le 1/\delta.
\end{equation}
The same estimate holds when $\theta$ is replaced by $-\theta$.

Next we estimate the norm of $Q_n^{\rm ee}$. On the $n^*$ nodes $i\alpha$
consider the matrix
\begin{equation}\label{eq:six-Z-base}
 Z_{ij}:=\frac{1}{2\sin\pi((i-j)\alpha)}\quad(i\ne j),\qquad Z_{ii}:=0 .
\end{equation}
Lemma~\ref{lem:six-cosecant} gives $\|Z\|\le1/(2\delta)$. To compare
$Q_n^{\rm ee}$ with $Z$, first note that
\begin{equation}\label{eq:six-sine-lower}
 |\sin(\pi k\alpha)|\ge2\|k\alpha\|_{\T}\ge2\delta
 \qquad(1\le |k|<n).
\end{equation}
The points $0,\alpha,\ldots,(n-1)\alpha$ are $\delta$-separated modulo $1$.
Ordering their nonzero representatives in $(-1/2,1/2]$ on each side of
zero, the $r$th closest on either side has distance at least $r\delta$.
Thus
\[
 \sum_{k=1}^{n-1}\frac{1}{\sin^2(\pi k\alpha)}
 \le\frac{2}{4\delta^2}\sum_{r=1}^\infty r^{-2}
 =\frac{\pi^2}{12\delta^2}<\frac{1}{\delta^2}.
\]
Since $|\sin\pi(k\alpha-\theta)-\sin(\pi k\alpha)|\le\pi|\theta|$,
\eqref{eq:six-sine-lower} and $|\theta|\le\delta/200$ give
$|\sin\pi(k\alpha-\theta)|\ge(9/10)|\sin(\pi k\alpha)|$. In particular,
for $i\ne j$ and $k=i-j$ we have
\[
 |Z_{ij}-(Q_n^{\rm ee})_{ij}|
 \le \frac{5\pi|\theta|}{9\sin^2(\pi k\alpha)}
 \le \frac{2|\theta|}{\sin^2(\pi k\alpha)}.
\]
Both diagonals vanish, and in a fixed row or column each $|k|$ occurs
at most twice. The preceding sum therefore gives
\[
 \|Z-Q_n^{\rm ee}\|_1,\ \|Z-Q_n^{\rm ee}\|_\infty
 \le \frac{4|\theta|}{\delta^2}\le\frac{1}{50\delta}.
\]
By \eqref{eq:row-column-bound} the same estimate holds for the operator
norm. It follows that
\begin{equation}\label{eq:six-q-reduced-bound}
 \|Q_n^{\rm ee}\|\le\|Z\|+\|Z-Q_n^{\rm ee}\|
 \le\frac{1}{2\delta}+\frac{1}{50\delta}\le \frac{2}{3\delta}.
\end{equation}
Finally $|b_0(\alpha;\theta)|\le2\pi|\theta|$. Combining
\eqref{eq:six-reduced-norms}, \eqref{eq:six-c-reduced-bound} and
\eqref{eq:six-q-reduced-bound} gives the two assertions.
\end{proof}

\subsubsection{Proof of Lemma~\ref{lem:six-criterion}}
\begin{proof}
By \eqref{eq:six-CQ-blocks}, the diagonal blocks of
$\Lambda_n(\alpha;\theta)-\mathbb J_n$ have norms
$|b_0(\alpha;\theta)|\,\|C_n(\alpha;\pm\theta)\|$.
The two off-diagonal blocks are $Q_n-I_n$ and its negative transpose,
so together they have norm $\|Q_n-I_n\|$. Lemma~\ref{lem:six-block-bounds}
now gives
\[
 \begin{aligned}
 \|\Lambda_n(\alpha;\theta)-\mathbb J_n\|
 &\le |b_0(\alpha;\theta)|\max_{\gamma=\pm\theta}\|C_n(\alpha;\gamma)\|
          +\|Q_n(\alpha;\theta)-I_n\|\\
 &\le 20|\theta|/\delta<1/5 .
 \end{aligned}
\]
Since $\mathbb J_n$ is orthogonal, the matrix
$A:=\mathbb J_n^T(\mathbb J_n-\Lambda_n(\alpha;\theta))$ has norm at most
$1/5$, and $\Lambda_n(\alpha;\theta)=\mathbb J_n(I_{2n}-A)$ is invertible.
Equation~\eqref{eq:six-lambda} then implies that $\Omega_n(\alpha;\theta)$
is invertible as well. By \eqref{eq:six-delta} and the Neumann series,
\[
 \Delta_n(\alpha;\theta)=\Lambda_n(\alpha;\theta)^{-1}
    =\left(\sum_{r=0}^\infty A^r\right)\mathbb J_n^T,\qquad
 \|\Delta_n(\alpha;\theta)\|\le \sum_{r=0}^\infty (1/5)^r=5/4 .
 \qedhere
\]
\end{proof}

}
\section{Proof of Lemma~\ref{lem:returnmatrix}}\label{sec:returnnorm}

\subsection{A uniform criterion for norm estimates}
At first we present the following uniform criterion for the norm estimates.

\begin{lemma}\label{lem:uniform-return}
There exists an absolute constant $\varrho>0$ such that for any $\alpha\in\mathbb I$,
if there exist integers $0<p<q$ satisfying $(p,q)=1$ and
\begin{equation}\label{eq:seven-approx}
 \left|\alpha-\frac pq\right|\le \frac{1}{A_0q^2},
\end{equation}
then for any $1\le N\le2q$,
\[
 |b_q(\alpha)|\,\|H^{[N]}(\alpha)\|\le\varrho .
\]
\end{lemma}

Lemma~\ref{lem:returnmatrix} is a direct consequence of this uniform criterion.

\begin{proof}[Proof of Lemma~\ref{lem:returnmatrix}]
{Here $0<p_k<q_k$ since $\alpha\in(0,1)$ and $k\ge20$.}
Since $k\in\mathcal I_\alpha$, we have $a_{k+1}\ge A_0$. Hence by
\eqref{eq:convergent-bounds},
\[
 \left|\alpha-\frac{p_k}{q_k}\right|
 \le\frac1{q_kq_{k+1}}\le\frac1{A_0q_k^2}.
\]
Since $m_k=q_k+n_k=q_k+\lfloor q_k/8\rfloor\le2q_k$, by
Lemma~\ref{lem:uniform-return} and \eqref{eq:Gamma-def},
\[
 \|\Gamma_k\|=\|b_{q_k}\Pi_{I_k}H^{[m_k]}\Pi_{m_k}\|
 \le |b_{q_k}|\,\|H^{[m_k]}\|\le\varrho .\qedhere
\]
\end{proof}

In the rest of this section, we prove Lemma~\ref{lem:uniform-return}. The plan is as follows.

At first, we decompose $H^{[N]}(\alpha)$ into a rational approximation part
$H^{[N]}(p/q)$, a resonant part $R^{[N]}$ and an error part $E^{[N]}$.
Then we estimate them separately. It is relatively easy to bound the norms
of $R^{[N]}$ and $E^{[N]}$. To bound the norm of $H^{[N]}(p/q)$,
we need to observe that it exhibits a Hankel structure, and then apply
some known Hankel truncation estimates.

In the rest of this section, we always write
\[
 \widetilde\alpha:=p/q.
\]

\subsection{The decomposition of \texorpdfstring{$H^{[N]}(\alpha)$}{H(N)}}
Naively, we just define the rational approximation part $H^{[N]}(\widetilde\alpha)$
by \eqref{eq:Hentries} with $\alpha$ replaced by $\widetilde\alpha$.
However by \eqref{eq:bcoeff}, $b_n(\widetilde\alpha)=0$ if $q\mid n$.
So we need to be more careful.

From now on, we always fix $\alpha\in\mathbb I$, $1\le N\le2q$.
Then we write
\[
 H:=H^{[N]}(\alpha).
\]
Recall from \eqref{eq:matrix-indices} that
\[
 \kappa_{ij}=\frac{i+j+1}{2}\quad\text{and}\quad \sigma_{ij}=\frac{i-j-1}{2}.
\]
Since $0\le i,j\le N-1\le2q-1$, we have
\begin{equation}\label{eq:seven-range}
 \tfrac12\le\kappa_{ij}<2q .
\end{equation}
{For $i>j$ with $i-j$ odd, $\kappa_{ij}$ is a positive integer,}
so $q\mid\kappa_{ij}$ if and only if $q=\kappa_{ij}$. Define the rational
approximation part $\widetilde H=H^{[N]}(\widetilde\alpha)$ as
\begin{equation}\label{eq:seven-rational}
 \widetilde H_{ij}=(H^{[N]}(\widetilde\alpha))_{ij}:=
 \begin{cases}
  (-1)^{\sigma_{ij}}/b_{\kappa_{ij}}(\widetilde\alpha),
       &i>j,\ i-j\text{ odd},\ \kappa_{ij}\ne q,\\
  0,&\text{otherwise}.
 \end{cases}
\end{equation}
Comparing \eqref{eq:Hentries} with \eqref{eq:seven-rational}, it is
reasonable to define the resonant part $R=R^{[N]}$ as
\begin{equation}\label{eq:seven-resonant}
 R_{ij}=(R^{[N]})_{ij}{{}:=}
 \begin{cases}
  (-1)^{\sigma_{ij}}/b_q(\alpha), &i>j,\ i-j\text{ odd},\ \kappa_{ij}=q,\\
  0,&\text{otherwise}.
 \end{cases}
\end{equation}
Finally, we simply define the error part $E=E^{[N]}$ as
\begin{equation}\label{eq:seven-decomposition}
 E=E^{[N]}:=H^{[N]}(\alpha)-H^{[N]}(\widetilde\alpha)-R^{[N]}=H-\widetilde H-R .
\end{equation}
Combining \eqref{eq:Hentries}, \eqref{eq:seven-rational} and \eqref{eq:seven-resonant},
\begin{equation}\label{eq:seven-errorentries}
 E_{ij}=(E^{[N]})_{ij}=
 \begin{cases}
  (-1)^{\sigma_{ij}}\big(b_{\kappa_{ij}}^{-1}(\alpha)-b_{\kappa_{ij}}^{-1}(\widetilde\alpha)\big),
      &i>j,\ i-j\text{ odd},\ \kappa_{ij}\ne q,\\
  0,&\text{otherwise}.
 \end{cases}
\end{equation}

\subsection{The estimates of \texorpdfstring{$\|R^{[N]}\|$ and $\|E^{[N]}\|$}{the resonant and error norms}}
The estimates of the norms of the resonant part and the error part are relatively easy.

\begin{lemma}\label{lem:resonance-error}
The following estimates hold:
\begin{equation}\label{eq:seven-two-bounds}
 |b_q(\alpha)|\,\|R^{[N]}\|\le1
 \qquad\text{and}\qquad
 |b_q(\alpha)|\,\|E^{[N]}\|\le1 .
\end{equation}
\end{lemma}

\begin{proof}
At first, we estimate $\|R^{[N]}\|$. By \eqref{eq:seven-resonant},
$R_{ij}\ne0$ only if $\kappa_{ij}=q$, or equivalently, $i+j+1=2q$.
This means that in each row and each column of $R^{[N]}$, there is at
most one nonzero element, which has {modulus} $1/|b_q(\alpha)|$.
Then the first estimate of \eqref{eq:seven-two-bounds} follows immediately.

Next, we estimate $\|E^{[N]}\|$. By \eqref{eq:bcoeff} and \eqref{eq:seven-approx},
for any integer $k\in[1,2q-1]$ we have
\[
 \begin{aligned}
 |b_k(\alpha)-b_k(\widetilde\alpha)|
 &=2|\sin(k\pi\alpha)-\sin(k\pi\widetilde\alpha)|
 \le4\left|\sin\frac{k\pi(\alpha-\widetilde\alpha)}2\right|\\
 &\le4q\pi|\alpha-\widetilde\alpha|\le\frac{4q\pi}{A_0q^2}\le\frac1{10q}.
 \end{aligned}
\]
For $k\in[1,2q-1]$ and $k\ne q$,
\[
 |b_k(\widetilde\alpha)|=2|\sin(kp\pi/q)|\ge2\sin(\pi/q)\ge4/q .
\]
So we get
\begin{equation}\label{eq:seven-denom}
 |b_k(\alpha)|\ge |b_k(\widetilde\alpha)|-\frac1{10q}
 \ge |b_k(\widetilde\alpha)|-\frac{|b_k(\widetilde\alpha)|}{40}
 \ge\frac{|b_k(\widetilde\alpha)|}2 .
\end{equation}
As a result, for $k\in[1,2q-1]$ and $k\ne q$, we have
\[
 |b_k^{-1}(\alpha)-b_k^{-1}(\widetilde\alpha)|
 =\frac{|b_k(\alpha)-b_k(\widetilde\alpha)|}{|b_k(\alpha)|\,|b_k(\widetilde\alpha)|}
 \le\frac{1}{5q|b_k(\widetilde\alpha)|^2}=\frac{\csc^2(k\pi\widetilde\alpha)}{20q}.
\]
It is known that (see {\cite[p.~99, Art.~37]{Jeffery1864}})
\[
 \sum_{k=1}^{q-1}\csc^2(k\pi/q)=\frac{q^2-1}{3}.
\]
Since $\widetilde\alpha=p/q$ and $(p,q)=1$, we also have
\[
 \sum_{k=1}^{q-1}\csc^2(k\pi\widetilde\alpha)
 =\sum_{k=q+1}^{2q-1}\csc^2(k\pi\widetilde\alpha)=\frac{q^2-1}{3}.
\]
Along any fixed row or column, the indices $\kappa_{ij}$ of the
nonzero entries are distinct. Now by \eqref{eq:seven-errorentries} and
\eqref{eq:seven-range} we have
\[
 \begin{aligned}
 \|E^{[N]}\|_1,\ \|E^{[N]}\|_\infty
 &\le\sum_{\substack{1\le k<2q\\ k\ne q}}
 |b_k^{-1}(\alpha)-b_k^{-1}(\widetilde\alpha)|\\
 &\le\frac1{20q}\sum_{\substack{1\le k<2q\\ k\ne q}}\csc^2(k\pi\widetilde\alpha)
 \le\frac1{20q}\cdot\frac{2q^2}{3}=\frac{q}{30} .
 \end{aligned}
\]
Then by \eqref{eq:row-column-bound}, $\|E^{[N]}\|\le q/30$.
On the other hand,
\begin{equation}\label{eq:seven-bq}
 |b_q(\alpha)|=2|\sin(q\pi\alpha)|=2|\sin(q\pi(\alpha-\widetilde\alpha))|
 \le2q\pi|\alpha-\widetilde\alpha|\le 1/q .
\end{equation}
So the second estimate in \eqref{eq:seven-two-bounds} holds.
\end{proof}

\subsection{The estimate of \texorpdfstring{$\|H^{[N]}(\widetilde\alpha)\|$}{the rational model norm} and proof of Lemma~\ref{lem:uniform-return}}
We have mentioned in the beginning of this section that the key to the norm estimates
of $H^{[N]}(\widetilde\alpha)$ is the hidden Hankel structure in it.
Let us first discuss the following concrete example:
\[
 H^{[9]}(1/5)=
 \begin{pmatrix}
 0&0&0&0&0&0&0&0&0\\
 c_1&0&0&0&0&0&0&0&0\\
 0&c_2&0&0&0&0&0&0&0\\
 -c_2&0&c_3&0&0&0&0&0&0\\
 0&-c_3&0&c_4&0&0&0&0&0\\
 c_3&0&-c_4&0&\mathbf{0}&0&0&0&0\\
 0&c_4&0&\mathbf{0}&0&c_6&0&0&0\\
 -c_4&0&\mathbf{0}&0&-c_6&0&c_7&0&0\\
 0&\mathbf{0}&0&c_6&0&-c_7&0&c_8&0
 \end{pmatrix},
\]
where $c_n=b_n^{-1}(1/5)$ for $1\le n\le9$, $n\ne5$, and $c_5=0$.
The bold zeros mark the resonant sites.

From this we already see the Hankel pattern modulo certain sign structure.
Now we do certain transformations such that the Hankel structure is more compact
in the sense that we gather the nonzero terms in two blocks. Let $P$ be the
permutation matrix related to the permutation $\sigma=(0,2,4,6,8,1,3,5,7)$.
Then we have
\[
 PH^{[9]}(1/5)P^T=\begin{pmatrix}0&\widehat H_{eo}\\\widehat H_{oe}&0\end{pmatrix},
\]
where the two matrices $\widehat H_{eo}$ and $\widehat H_{oe}$ are given by
\[
 \widehat H_{eo}=\begin{pmatrix}
 0&0&0&0\\c_2&0&0&0\\-c_3&c_4&0&0\\c_4&-c_5&c_6&0\\-c_5&c_6&-c_7&c_8
 \end{pmatrix}
 \quad\text{and}\quad
 \widehat H_{oe}=\begin{pmatrix}
 c_1&0&0&0&0\\-c_2&c_3&0&0&0\\c_3&-c_4&c_5&0&0\\-c_4&c_5&-c_6&c_7&0
 \end{pmatrix}.
\]
We can even make all the {prefactors of the $c_n$} positive. Let
\[
 \begin{aligned}
 D_{1l}&:=\operatorname{diag}(1,1,-1,1,-1),&D_{1r}&:=\operatorname{diag}(1,-1,1,-1),\\
 D_{2l}&:=\operatorname{diag}(1,-1,1,-1),&D_{2r}&:=\operatorname{diag}(1,-1,1,-1,1).
 \end{aligned}
\]
Then
\[
 H_{eo}:=D_{1l}\widehat H_{eo}D_{1r}=
 \begin{pmatrix}
 0&0&0&0\\c_2&0&0&0\\c_3&c_4&0&0\\c_4&c_5&c_6&0\\c_5&c_6&c_7&c_8
 \end{pmatrix},
 \qquad
 H_{oe}:=D_{2l}\widehat H_{oe}D_{2r}=
 \begin{pmatrix}
 c_1&0&0&0&0\\c_2&c_3&0&0&0\\c_3&c_4&c_5&0&0\\c_4&c_5&c_6&c_7&0
 \end{pmatrix}.
\]
It is seen that
\[
 \|H^{[9]}(1/5)\|=\max\{\|\widehat H_{eo}\|,\|\widehat H_{oe}\|\}
   =\max\{\|H_{eo}\|,\|H_{oe}\|\}.
\]
Now if we define a Hankel matrix
\[
 C:=\begin{pmatrix}
 c_1&c_2&c_3&c_4&c_5\\c_2&c_3&c_4&c_5&c_6\\
 c_3&c_4&c_5&c_6&c_7\\c_4&c_5&c_6&c_7&c_8\\c_5&c_6&c_7&c_8&c_9
 \end{pmatrix},
\]
then we see that both $H_{eo}$ and $H_{oe}$ are truncations of $C$.
More precisely,
{
\[
 H_{eo}=(\chi_{i\ge j+1}c_{i+j+1})_{0\le i<5,\,0\le j<4},
 \qquad
 H_{oe}=(\chi_{i\ge j}c_{i+j+1})_{0\le i<4,\,0\le j<5}.
\]}
This makes it possible to apply some known results \cite{BB} to estimate the norms
of $H_{eo}$ and $H_{oe}$. To do this we need to relate
the sequence $\{c_n\}$ to the Fourier coefficients of a certain continuous function
{$\psi_s$} on $\mathbb T$. Notice that $\{c_n\}$ is a nontrivial periodic sequence,
so it cannot be the Fourier coefficients of any continuous function on $\mathbb T$.
We must multiply the $c_n$ by a decaying factor. This motivates the definition in Section~\ref{sec:seven-symbol}.

The general case is exactly the same as we describe above. Now we give the full detail.

\subsubsection{The Hankel structure of \texorpdfstring{$H^{[N]}(\widetilde\alpha)$}{the rational model}}
Define a sequence $\{c_n:n\in\mathbb Z\}$ as
\begin{equation}\label{eq:seven-c}
 c_n:=\begin{cases}b_n^{-1}(\widetilde\alpha),&q\nmid n,\\0,&q\mid n.\end{cases}
\end{equation}
Since $b_n(\widetilde\alpha)=2\sin(np\pi/q)$, the sequence
$(c_n)_{n\in\mathbb Z}$ is $2q$-periodic.
Define an infinite Hankel matrix $C(\widetilde\alpha)$ as
\begin{equation}\label{eq:seven-C}
 C(\widetilde\alpha)_{ij}:=c_{i+j+1},\qquad i,j\in\mathbb Z^+ .
\end{equation}
For any matrix $A=(a_{ij})_{i,j\ge0}$ and $l\ge0$, we define the $l$-th
lower triangular truncation of $A$ as
\begin{equation}\label{eq:seven-trunc}
 (T_lA)_{ij}:=\chi_{\{i\ge j+l\}}a_{ij}.
\end{equation}
For any $M,N\ge0$, we define
\[
 A[M,N]:=(a_{ij})_{0\le i<M,\,0\le j<N}.
\]
Now we have the following structure property:
\begin{lemma}\label{lem:rational-structure}
For any $1\le N\le2q$, there exist unitary matrices $U$ and $V$ such that
\[
 H^{[N]}(\widetilde\alpha)=
 U\begin{pmatrix}
 0&H_{eo}^{[N]}\\
 H_{oe}^{[N]}&0
 \end{pmatrix}V,
\]
where
\[
 \begin{aligned}
 H_{eo}^{[N]}&:=(T_1C(\widetilde\alpha))[\lceil N/2\rceil,\lfloor N/2\rfloor],\\
 H_{oe}^{[N]}&:=(T_0C(\widetilde\alpha))[\lfloor N/2\rfloor,\lceil N/2\rceil].
 \end{aligned}
\]
\end{lemma}
{
\begin{proof}
Order the even coordinates first and the odd coordinates second. The diagonal
blocks vanish by \eqref{eq:seven-rational}. In the odd-row, even-column block
put $i=2m+1$ and $j=2n$. In the even-row, odd-column block put $i=2m$ and $j=2n+1$.
The respective entries of the two blocks before the sign changes are
\[
 (\widehat H_{oe}^{[N]})_{mn}=\chi_{\{m\ge n\}}(-1)^{m-n}c_{m+n+1},\qquad
 (\widehat H_{eo}^{[N]})_{mn}=\chi_{\{m\ge n+1\}}(-1)^{m-n-1}c_{m+n+1}.
\]
Multiply odd row $2m+1$ by $(-1)^m$, even row $2m$ by $(-1)^{m+1}$,
and both columns $2n$ and $2n+1$ by $(-1)^n$ whenever these indices occur.
This removes both signs and gives $H_{oe}^{[N]}$ and $H_{eo}^{[N]}$, respectively. There are $\lceil N/2\rceil$ even coordinates and
$\lfloor N/2\rfloor$ odd coordinates. Undoing the signs and permutations
therefore gives the required unitary matrices.
\end{proof}}

\subsubsection{The auxiliary function \texorpdfstring{$\psi_s$}{psi(s)} and the associated Hankel matrix}\label{sec:seven-symbol}
Fix any $0<s<1$ and define a periodic function $\psi_s$ on $\mathbb R$ as
\begin{equation}\label{eq:seven-symbol}
 \psi_s(x):=\sum_{n\in\mathbb Z}s^{|n|}c_{n+1}\mathrm e^{2\pi n x\mathrm i}.
\end{equation}
Since $0<s<1$, the above series converges both absolutely and uniformly.
Consequently $\psi_s$ is continuous and $1$-periodic. It has the following
basic property.
\begin{lemma}\label{lem:symbol-bound}
The function $\psi_s$ satisfies
\[
 \|\psi_s\|_\infty\le\frac q2\frac{1+s^q}{1-s^q},
 \qquad \widehat{\psi_s}(n)=s^{|n|}c_{n+1}.
\]
\end{lemma}

{
\begin{proof}
Termwise integration in \eqref{eq:seven-symbol} gives the Fourier coefficient
with the convention $\widehat{\psi_s}(n)=\int_0^1\psi_s(x)\mathrm e^{-2\pi\mathrm i n x}\,dx$.
To estimate the supremum norm, put $g_r:=\mathrm e^{\pi\mathrm i p r/q}c_r$.
The relation $c_{r+q}=(-1)^p c_r$ makes $g_r$ periodic of period $q$.
Write its finite Fourier expansion as
\[
 g_r=\sum_{\ell=0}^{q-1}\omega_\ell\,\mathrm e^{2\pi\mathrm i\ell r/q}.
\]
By \eqref{eq:seven-c}, $g_r(1-\mathrm e^{-2\pi\mathrm i p r/q})
=\mathrm i(1-\chi_{\{q\mid r\}})$ for every $r$. Comparing Fourier coefficients
and evaluating the expansion at $r=0$ gives
\[
 \omega_\ell-\omega_{\ell+p}=\mathrm i(\chi_{\{\ell=0\}}-1/q),
 \qquad \sum_{\ell=0}^{q-1}\omega_\ell=0,
\]
where the subscripts are taken modulo $q$. Since multiplication by $p$ permutes the
residues, the first relation gives $\omega_{jp}=\omega_0-\mathrm i+\mathrm i j/q$
for $1\le j<q$. The second now gives
\[
 \omega_0=\mathrm i\,\frac{q-1}{2q},\qquad
 \omega_{jp}=\mathrm i\,\frac{2j-q-1}{2q}\quad(1\le j<q).
\]
In particular $|\omega_\ell|<1/2$.

Put $\beta_\ell:=\ell/q-p/(2q)$. Then $c_r=\sum_\ell\omega_\ell
\mathrm e^{2\pi\mathrm i r\beta_\ell}$. Substituting in \eqref{eq:seven-symbol} yields
\[
 \psi_s(x)=\sum_{\ell=0}^{q-1}\omega_\ell\mathrm e^{2\pi\mathrm i\beta_\ell}
 P_s(x+\beta_\ell),\qquad
 P_s(x):=\sum_{n\in\mathbb Z}s^{|n|}\mathrm e^{2\pi\mathrm i nx}
 =\frac{1-s^2}{1-2s\cos(2\pi x)+s^2}.
\]
In the sum of these Poisson kernels, the sum over $\ell$ vanishes at all
Fourier indices except the multiples of $q$. Hence
\[
 |\psi_s(x)|\le\frac12\sum_{\ell=0}^{q-1}P_s(x+\beta_\ell)
 =\frac q2 P_{s^q}(qx-p/2)
 \le\frac q2\,\frac{1+s^q}{1-s^q}.\qedhere
\]
\end{proof}}

\subsubsection{Hankel truncation estimate and its consequence}
We now use the truncation theorem of Bonami and Bruna to estimate
$\|H^{[N]}(\widetilde\alpha)\|$.
For $\psi\in L^\infty(\mathbb T)$, define the Hankel matrix $\Gamma_\psi$ by
\begin{equation}\label{eq:seven-Gamma}
 (\Gamma_\psi)_{ij}:=\widehat\psi(i+j),\qquad i,j\ge0.
\end{equation}

\begin{proposition}[{\cite[Theorem~1(a) and (2) on p.~238]{BB}}]\label{prop:seven-trunc}
The matrix $\Gamma_\psi$ defines a bounded operator on $\ell^2(\mathbb Z^+)$
with $\|\Gamma_\psi\|\le\|\psi\|_\infty$. There exists an absolute constant
$\vartheta>0$ such that for every $l\ge0$,
\[
 \|T_l\Gamma_\psi\|\le\vartheta\|\psi\|_\infty .
\]
\end{proposition}

\begin{proof}
For a finitely supported sequence $x$, put
$f(t)=\sum_{j\ge0}x_j\mathrm e^{-2\pi\mathrm i jt}$. Then
$(\Gamma_\psi x)_i=\widehat{\psi f}(i)$ for $i\ge0$, and Parseval's
identity gives $\|\Gamma_\psi x\|\le\|\psi\|_\infty\|x\|$.
For the truncation, equation~(2) of \cite{BB} with slope $1$ and intercept $l$
gives
\[
 \sum_{i\ge0}\left|\sum_{0\le j\le i-l}\widehat\psi(i+j)x_j\right|^2
 \le C(1)\|\psi\|_\infty^2\sum_{j\ge0}|x_j|^2 .
\]
The constant is independent of $l$. Taking square roots and using density
gives the asserted bound with $\vartheta=\sqrt{C(1)}$.
\end{proof}

Now put $\Gamma_s:=\Gamma_{\psi_s}$. By Lemma~\ref{lem:symbol-bound},
\begin{equation}\label{eq:seven-Gammas}
 (\Gamma_s)_{ij}=\widehat{\psi_s}(i+j)=s^{i+j}c_{i+j+1},
 \qquad i,j\ge0 .
\end{equation}
Writing $D_s:=\operatorname{diag}(1,s,s^2,\ldots)$, we obtain the
entrywise identity
\begin{equation}\label{eq:seven-weight}
 \Gamma_s=D_sC(\widetilde\alpha)D_s .
\end{equation}
Through this identity we can estimate both blocks
$H_{eo}^{[N]}$ and $H_{oe}^{[N]}$.

\begin{lemma}\label{lem:seven-rational-bound}
For any $1\le N\le2q$,
\[
 |b_q(\alpha)|\,\|H^{[N]}(\widetilde\alpha)\|\le12\vartheta .
\]
\end{lemma}
\begin{proof}
When $N=1$ the matrix is zero. For $N\ge2$, put $m:=\lceil N/2\rceil$.
The two blocks in Lemma~\ref{lem:rational-structure} act on orthogonal
input coordinates and have orthogonal images, so
\begin{equation}\label{eq:seven-blocknorm}
 \|H^{[N]}(\widetilde\alpha)\|
 =\max\big\{\|H_{eo}^{[N]}\|,\|H_{oe}^{[N]}\|\big\}.
\end{equation}
Both blocks are leading rectangles of size $a\times b$ with $a,b\le m$.
For $l=1$ in the even-row block and $l=0$ in the odd-row block,
\eqref{eq:seven-weight} gives
\[
 (T_l C(\widetilde\alpha))[a,b]=
 \operatorname{diag}(s^{-i})_{0\le i<a}\,
 (T_l\Gamma_s)[a,b]\,
 \operatorname{diag}(s^{-j})_{0\le j<b}.
\]
Extending the input by zeros and retaining the first $a$ output coordinates
gives $\|(T_l\Gamma_s)[a,b]\|\le\|T_l\Gamma_s\|$. The product of the two
diagonal norms is $s^{-(a+b-2)}\le s^{-2m}$. By Proposition~\ref{prop:seven-trunc}
and Lemma~\ref{lem:symbol-bound},
\[
 \|(T_l C(\widetilde\alpha))[a,b]\|
 \le s^{-2m}\vartheta\|\psi_s\|_\infty
 \le \frac{\vartheta q(1+s^q)}{2s^{2m}(1-s^q)} .
\]
Take $s=2^{-1/q}$. Then $s^q=1/2$ and $s^{-2m}=2^{2m/q}\le4$ since
$m\le q$. Thus each block has norm at most $6\vartheta q$, in particular
at most $12\vartheta q$. Now \eqref{eq:seven-blocknorm} and
$|b_q(\alpha)|\le 1/q$ from \eqref{eq:seven-bq} give the result.
\end{proof}

\subsubsection{Proof of Lemma~\ref{lem:uniform-return}}
\begin{proof}
By \eqref{eq:seven-decomposition}, Lemma~\ref{lem:resonance-error}
and Lemma~\ref{lem:seven-rational-bound},
\[
 \begin{aligned}
 |b_q(\alpha)|\,\|H^{[N]}(\alpha)\|
 &\le |b_q(\alpha)|\big(\|R^{[N]}\|+\|E^{[N]}\|
                      +\|H^{[N]}(\widetilde\alpha)\|\big)\\
 &\le 1+1+12\vartheta .
 \end{aligned}
\]
Taking $\varrho:=12\vartheta+2$ proves the lemma.
\end{proof}
\section{Proof of Lemma~\ref{lem:six-entries}}\label{sec:six-tech}
By \eqref{eq:six-CQ}, to obtain the formulas for $C_n(\alpha;\theta)$ and
$Q_n(\alpha;\theta)$, the key step is to find the formulas for
$U(\alpha;\theta)$ and $V(\alpha;\theta)$.

\subsection{The formulas for \texorpdfstring{$U(\alpha;\theta)$ and $V(\alpha;\theta)$}{U and V}}\label{subsec:eight-connections}
Fix $n\in\N$. Assume that the phase $\theta$ satisfies
\begin{equation}\label{eq:eight-one-phase}
 b_k(\alpha;\theta)\ne0\qquad(1\le k<n).
\end{equation}
The polynomial expansions in the proof of Lemma~\ref{lem:six-connection}
define both $n\times n$ connection matrices at this phase. Put
\begin{equation}\label{eq:eight-varpi}
 \varpi_0(\alpha;\theta):=1,\qquad
 \varpi_k(\alpha;\theta):=\prod_{s=1}^k\frac{b_s(\alpha)}{b_s(\alpha;\theta)}
 \qquad(1\le k<n).
\end{equation}
\begin{lemma}\label{lem:eight-coefficients}
Assume \eqref{eq:eight-one-phase}. Then for $0\le i,j<n$,
\begin{equation}\label{eq:U}
 (U(\alpha;\theta))_{ij}=
 \begin{cases}
 0,&i<j\ \text{or }i-j\text{ odd},\\
 \varpi_i(\alpha;\theta),&i=j,\\[2pt]
 \displaystyle(-1)^h
 \frac{b_0(\alpha;\theta)\varpi_{i-h}(\alpha;\theta)}
      {b_h(\alpha;\theta)\varpi_h(\alpha;\theta)},&i-j=2h>0.
 \end{cases}
\end{equation}
For the second matrix,
\begin{equation}\label{eq:Vshift}
 (V(\alpha;\theta))_{ij}=
 \begin{cases}
 0,&i\le j,\\[2pt]
 \displaystyle\frac{(U(\alpha;\theta+\alpha))_{i-1,j}}{b_1(\alpha;\theta)},
      & i>j .
 \end{cases}
\end{equation}
The matrix at $\theta+\alpha$ in the second line has size $n-1$.
\end{lemma}
\begin{proof}
For $n=1$ the matrices are $(1)$ and $(0)$. Assume $n\ge2$.
Write $b_j=b_j(\alpha)$. For a phase $\gamma$, let
$\varpi_j(\gamma):=\varpi_j(\alpha;\gamma)$ and set
\[
 s_j(\gamma):=b_j(\alpha;\gamma)=2\sin\!\bigl(\pi(j\alpha+\gamma)\bigr).
\]
Denote the finite connection matrices by $U_\gamma,V_\gamma$ at the size
indicated by their indices. The parity of the recurrences and the polynomial
expansions give
\begin{equation}\label{eq:connectionparity}
 (U_\gamma)_{ij}=0\quad(i-j\text{ odd}),\qquad
 (V_\gamma)_{ij}=0\quad(i-j\text{ even}).
\end{equation}
For the calculation at $\gamma=\theta$, write $s_j=s_j(\gamma)$, $\varpi_j=\varpi_j(\gamma)$ and $u_i^\gamma=u_i(z,\alpha;\gamma)$, and likewise for $v_i^\gamma$. Comparing leading coefficients in
$u_i^\gamma=\sum_{j\le i}(U_\gamma)_{ij}u_j$ gives
$(U_\gamma)_{ii}=(b_1\cdots b_i)/(s_1\cdots s_i)=\varpi_i$.

Assume first $s_0\ne0$, which permits the divisions by $s_0$ below. Let $U$ be the lower triangular matrix with
diagonal $\varpi_i$, the entries \eqref{eq:U} for $i-j=2h>0$, and zeros for
odd $i-j$; then the diagonal is also given by \eqref{eq:U} with $h=0$,
read as $s_0/s_0=1$. Put $\mathsf R_i(z)=\sum_{j\le i}U_{ij}u_j(z)$. Then
$\mathsf R_0=1$ and $\mathsf R_1=\varpi_1u_1=z/s_1$, the initial values of
$u^\gamma$. Using $zu_l=b_lu_{l-1}+b_{l+1}u_{l+1}$, the coefficient of
$u_j$ in $z\mathsf R_i$ is $b_jU_{i,j-1}+b_{j+1}U_{i,j+1}$, and in
$s_i\mathsf R_{i-1}+s_{i+1}\mathsf R_{i+1}$ it is
$s_iU_{i-1,j}+s_{i+1}U_{i+1,j}$. So the recurrence for $\mathsf R$ reduces to
\begin{equation}\label{eq:eight-connection-recurrence}
 b_jU_{i,j-1}+b_{j+1}U_{i,j+1}=s_iU_{i-1,j}+s_{i+1}U_{i+1,j}
 \qquad(1\le i<n-1,\ j\ge0),
\end{equation}
with entries outside the index range read as $0$. If $i-j$ is even all four
terms vanish by parity, and they also vanish if $j>i+1$. If $j=i+1$, the
identity reads $b_{i+1}\varpi_i=s_{i+1}\varpi_{i+1}$, which holds by
definition. It remains to treat $i-j=2m+1$ with $m\ge0$. Put $\ell=i-m$,
so that $i=m+\ell$ and $j=\ell-m-1\ge0$. By \eqref{eq:U},
\[
\begin{aligned}
 U_{i,j-1}&=\frac{(-1)^{m+1}s_0\varpi_{\ell-1}}{s_{m+1}\varpi_{m+1}},&
 U_{i,j+1}&=\frac{(-1)^ms_0\varpi_\ell}{s_m\varpi_m},\\
 U_{i-1,j}&=\frac{(-1)^ms_0\varpi_{\ell-1}}{s_m\varpi_m},&
 U_{i+1,j}&=\frac{(-1)^{m+1}s_0\varpi_\ell}{s_{m+1}\varpi_{m+1}},
\end{aligned}
\]
the first being needed only for $j\ge1$. Dividing the difference of the
two sides of \eqref{eq:eight-connection-recurrence} by
$(-1)^ms_0\varpi_{\ell-1}/(s_ms_{m+1}\varpi_{m+1}s_\ell)$ and using
$\varpi_{t+1}/\varpi_t=b_{t+1}/s_{t+1}$, we obtain the expression
\begin{equation}\label{eq:eight-triple}
 b_\ell\bigl[b_{\ell-m}b_{m+1}+s_{m+\ell+1}s_m\bigr]
 -s_\ell\bigl[b_{\ell-m-1}s_m+s_{m+\ell}b_{m+1}\bigr].
\end{equation}
For $j=0$ the omitted term carries the factor $b_{\ell-m-1}=b_0=0$, so
\eqref{eq:eight-triple} applies as well. By the product-to-sum formula,
\[
\begin{aligned}
 b_{\ell-m}b_{m+1}+s_{m+\ell+1}s_m
 &=2\bigl[\cos\bigl(\pi(\ell-2m-1)\alpha\bigr)
   -\cos\bigl(\pi((\ell+2m+1)\alpha+2\gamma)\bigr)\bigr]\\
 &=s_\ell s_{2m+1},\\
 b_{\ell-m-1}s_m+s_{m+\ell}b_{m+1}
 &=2\bigl[\cos\bigl(\pi((\ell-2m-1)\alpha-\gamma)\bigr)
   -\cos\bigl(\pi((\ell+2m+1)\alpha+\gamma)\bigr)\bigr]\\
 &=b_\ell s_{2m+1},
\end{aligned}
\]
so \eqref{eq:eight-triple} vanishes. Hence $\mathsf R$ satisfies the
phase-$\gamma$ recurrence with the initial values of $u^\gamma$, and
$\mathsf R_i=u_i^\gamma$ for $i<n$. This proves \eqref{eq:U} when
$s_0\ne0$. If $s_0=0$, both sides of \eqref{eq:U} are continuous in the
phase on the open set where $s_1,\ldots,s_{n-1}$ are nonzero (the left side
because $u^\gamma_i$ depends continuously on $\gamma$ there), and the
zeros of $s_0$ are isolated; so \eqref{eq:U} extends.

For \eqref{eq:Vshift}, note that
$s_j(\gamma+\alpha)=s_{j+1}(\gamma)$. Define
\[
 y_0=0,\qquad
 y_i=s_1(\gamma)^{-1}u^{\gamma+\alpha}_{i-1}(z)
 \quad(1\le i<n).
\]
Then $y_0=0$ and $y_1=1/s_1(\gamma)$, the initial values of $v^\gamma$.
For $n\ge3$, we have $y_2=z/(s_1s_2)$, which verifies the recurrence at
$i=1$. For $2\le i<n-1$, the shifted recurrence gives
\[
 s_{i+1}(\gamma)y_{i+1}=zy_i-s_i(\gamma)y_{i-1}.
\]
Thus $y$ satisfies the phase-$\gamma$ recurrence throughout the finite
block, and uniqueness gives $y_i=v_i^\gamma$ for $0\le i<n$.
When $n=2$, the initial values already give this conclusion.
Expanding $u^{\gamma+\alpha}_{i-1}$ in the basis $(u_j)$ proves
\eqref{eq:Vshift}.
\end{proof}

\subsection{Proof of Lemma~\ref{lem:six-entries}}\label{subsec:eight-entries}
\begin{proof}
Let $n$ and $\theta$ satisfy the assumptions of Lemma~\ref{lem:six-entries}.
For $n=1$, $C_1(\alpha;\theta)=0$ and $Q_1(\alpha;\theta)=I_1$, so both
entry formulas hold. Suppose $n\ge2$. We use $b_j$, $s_j(\gamma)$ and
$\varpi_j(\gamma)$ as in Section~\ref{subsec:eight-connections} and write
$U_\gamma=U^{[n]}(\alpha;\gamma)$, $V_\gamma=V^{[n]}(\alpha;\gamma)$.
Here $\gamma=\theta$ or $-\theta$ unless another phase is specified.
The nonvanishing assumption ensures that both finite expansions in
Lemma~\ref{lem:six-connection} exist. Set
\[
 C_\gamma:=U_\gamma V_\gamma^T-V_\gamma U_\gamma^T,\qquad
 Q_\gamma:=U_{-\gamma}U_\gamma^T-s_0(\gamma)^2V_{-\gamma}V_\gamma^T.
\]
Thus, in the notation of \eqref{eq:six-CQ},
\[
 C_\gamma=C_n(\alpha;\gamma),\qquad
 Q_\gamma=Q_n(\alpha;\gamma)^T=Q_n(\alpha;-\gamma).
\]
We will obtain the following entry formulas:
\begin{equation}\label{eq:bilinearentries}
\begin{aligned}
 (C_\gamma)_{ij}
 &=\begin{cases}
 \dfrac{(-1)^{(i-j+1)/2}}{s_{(i+j+1)/2}(\gamma)},&i-j\ \text{odd},\\
 0,&i-j\ \text{even},
 \end{cases}\\
 (Q_\gamma)_{ij}
 &=\begin{cases}
 1,&i=j,\\
 \dfrac{(-1)^{h+1}s_0(\gamma)}{s_h(-\gamma)},&i-j=2h\ne0,\\
 0,&i-j\ \text{odd}.
 \end{cases}
\end{aligned}
\end{equation}
The entries exhibit cancellation between products of connection matrices.
We give the calculation, since bounds on the individual products would not
by themselves provide a bound independent of $n$.

\noindent\emph{The entries of $C_\gamma$.}
Assume for the moment that $s_0(\gamma)\ne0$. Combining \eqref{eq:U} and \eqref{eq:Vshift}
(and $s_0(\gamma+\alpha)=s_1(\gamma)$) gives, for $t,r\ge0$ and indices
below $n$,
\begin{equation}\label{eq:connectionproducts}
 (U_\gamma)_{r+2t,r}=(-1)^ts_0(\gamma)
 \frac{\prod_{\nu=t+1}^{t+r}b_\nu}{\prod_{\nu=t}^{t+r}s_\nu(\gamma)},\qquad
 (V_\gamma)_{r+2t+1,r}=(-1)^t
 \frac{\prod_{\nu=t+1}^{t+r}b_\nu}{\prod_{\nu=t+1}^{t+r+1}s_\nu(\gamma)},
\end{equation}
with empty products equal to $1$; for $t=0$ the first formula is
$(U_\gamma)_{rr}=\varpi_r(\gamma)$. All other entries on or below the
diagonal vanish by \eqref{eq:connectionparity}. For integers $a,a'$ the
product-to-sum formula gives
\begin{equation}\label{eq:samephasesine}
 s_a(\gamma)s_{a'}(\gamma)-b_ab_{a'}
 =4\sin(\pi\gamma)\sin\bigl(\pi((a+a')\alpha+\gamma)\bigr)
 =s_0(\gamma)s_{a+a'}(\gamma).
\end{equation}

By \eqref{eq:connectionparity}, $(C_\gamma)_{ij}=0$ for even $i-j$, and
$C_\gamma^T=-C_\gamma$. Let $i=j+2h+1$ with $h\ge0$, and put
$\ell=j+h+1=(i+j+1)/2\le n-1$. Order the indices $0,\ldots,j$ as
\begin{equation}\label{eq:sumorder}
 \rho_r=\frac{j-r}2\ \ (j-r\text{ even}),\qquad
 \rho_r=\frac{j+r+1}2\ \ (j-r\text{ odd}),\qquad 0\le r\le j .
\end{equation}
Put $t_r=\lfloor(j-r)/2\rfloor$ and
$\mathcal A_r=\{\rho_0,\ldots,\rho_{r-1}\}$, with $\mathcal A_0=\varnothing$.
For $0\le r\le j$, the definition gives
\[
 (\rho_r,t_{r+1})=
 \begin{cases}
 (t_r,t_r-1),&j-r\text{ even},\\
 (t_r+r+1,t_r),&j-r\text{ odd}.
 \end{cases}
\]
Each step adds one endpoint to the preceding interval. Induction yields
\[
 \mathcal A_r=\{t_r+1,\ldots,t_r+r\}\qquad(0\le r\le j+1).
\]
In particular, $t_{j+1}=-1$ and $\mathcal A_{j+1}=\{0,\ldots,j\}$, so
$(\rho_0,\ldots,\rho_j)$ is a permutation of these indices.
The definition also gives $\rho_j=0$. Define
\[
 \mathsf P_r=\prod_{\nu=0}^{r-1}\frac{b_{\rho_\nu}b_{\ell-\rho_\nu}}
 {s_{\rho_\nu}(\gamma)s_{\ell-\rho_\nu}(\gamma)}\qquad(0\le r\le j+1),
\]
so that $\mathsf P_0=1$ and $\mathsf P_{j+1}=0$, since its last factor
contains $b_{\rho_j}=b_0=0$. The temporary assumption $s_0(\gamma)\ne0$
makes that factor well defined; it also justifies the diagonal case
$t=0$ in \eqref{eq:connectionproducts}.
The $(i,j)$ entry of $C_\gamma$ is
$\sum_{r\le j}\bigl[(U_\gamma)_{ir}(V_\gamma)_{jr}-(V_\gamma)_{ir}(U_\gamma)_{jr}\bigr]$.
For the following two calculations write $s_\nu=s_\nu(\gamma)$.
If $j-r=2t$, only the second product is nonzero, and
\eqref{eq:connectionproducts} gives
\[
 -(V_\gamma)_{ir}(U_\gamma)_{jr}
 =(-1)^{h+1}s_0
 \frac{\displaystyle\prod_{\nu=t+1}^{t+r}b_\nu
       \prod_{\nu=h+t+1}^{h+t+r}b_\nu}
      {\displaystyle\prod_{\nu=t}^{t+r}s_\nu
       \prod_{\nu=h+t+1}^{h+t+r+1}s_\nu}.
\]
Here $\rho_r=t$, $\ell-\rho_r=h+t+r+1$, and
$\{\ell-\rho_\nu:0\le \nu<r\}=\{h+t+1,\ldots,h+t+r\}$.
If $j-r=2t+1$, only the first product is nonzero, and
\[
 (U_\gamma)_{ir}(V_\gamma)_{jr}
 =(-1)^{h+1}s_0
 \frac{\displaystyle\prod_{\nu=t+1}^{t+r}b_\nu
       \prod_{\nu=h+t+2}^{h+t+r+1}b_\nu}
      {\displaystyle\prod_{\nu=t+1}^{t+r+1}s_\nu
       \prod_{\nu=h+t+1}^{h+t+r+1}s_\nu}.
\]
In this case $\rho_r=t+r+1$, $\ell-\rho_r=h+t+1$, and
$\{\ell-\rho_\nu:0\le \nu<r\}=\{h+t+2,\ldots,h+t+r+1\}$.
In both cases the first interval of numerator indices is
$\{\rho_\nu:0\le \nu<r\}=\{t+1,\ldots,t+r\}$.
Comparing the products with the definition of $\mathsf P_r$ therefore yields
\begin{equation}\label{eq:Ctelescoping}
 (U_\gamma)_{ir}(V_\gamma)_{jr}-(V_\gamma)_{ir}(U_\gamma)_{jr}
 =\frac{(-1)^{h+1}s_0(\gamma)\mathsf P_r}{s_{\rho_r}(\gamma)s_{\ell-\rho_r}(\gamma)}
 =\frac{(-1)^{h+1}}{s_\ell(\gamma)}\bigl(\mathsf P_r-\mathsf P_{r+1}\bigr),
\end{equation}
where the last step is \eqref{eq:samephasesine} with $a=\rho_r$,
$a'=\ell-\rho_r$. Summing over $r=0,\ldots,j$ telescopes to
$(C_\gamma)_{ij}=(-1)^{h+1}/s_\ell(\gamma)=(-1)^{(i-j+1)/2}/s_{(i+j+1)/2}(\gamma)$.
For odd $i-j$, $(-1)^{(j-i+1)/2}=-(-1)^{(i-j+1)/2}$, so the formula is
antisymmetric in $(i,j)$, and skew-symmetry gives the upper triangle.

\medskip\noindent\emph{The entries of $Q_\gamma$.}
Continue to assume $s_0(\gamma)\ne0$. By parity,
$(Q_\gamma)_{ij}=0$ for odd $i-j$. Let $i=j+2h$ with $h\ge0$, use the
order \eqref{eq:sumorder}, and now put
\[
 \widetilde{\mathsf P}_r=\prod_{\nu=0}^{r-1}\frac{b_{h+\rho_\nu}b_{\rho_\nu}}
 {s_{h+\rho_\nu}(-\gamma)s_{\rho_\nu}(\gamma)}\qquad(0\le r\le j+1),
\]
so again $\widetilde{\mathsf P}_0=1$ and $\widetilde{\mathsf P}_{j+1}=0$; here
$h+\rho_\nu\le(i+j)/2\le n-1$. When $h=0$, the last factor also contains
$s_0(-\gamma)=-s_0(\gamma)$ in its denominator, which is nonzero under
the temporary assumption.
The analogue of \eqref{eq:samephasesine} for opposite phases is
\begin{equation}\label{eq:oppositephasesine}
 s_a(-\gamma)s_{a'}(\gamma)-b_ab_{a'}
 =4\sin(\pi\gamma)\sin\bigl(\pi((a-a')\alpha-\gamma)\bigr)
 =s_0(\gamma)s_{a-a'}(-\gamma).
\end{equation}
We again compute the two parity cases separately. If $j-r=2t$,
only the product of the $U$ entries contributes. Since
$s_0(-\gamma)=-s_0(\gamma)$, \eqref{eq:connectionproducts} gives
\[
 (U_{-\gamma})_{ir}(U_\gamma)_{jr}
 =(-1)^{h+1}s_0(\gamma)^2
 \frac{\displaystyle\prod_{\nu=h+t+1}^{h+t+r}b_\nu
       \prod_{\nu=t+1}^{t+r}b_\nu}
      {\displaystyle\prod_{\nu=h+t}^{h+t+r}s_\nu(-\gamma)
       \prod_{\nu=t}^{t+r}s_\nu(\gamma)}.
\]
Here $\rho_r=t$. If $j-r=2t+1$, only the product of the $V$ entries
contributes, and
\[
 -s_0(\gamma)^2(V_{-\gamma})_{ir}(V_\gamma)_{jr}
 =(-1)^{h+1}s_0(\gamma)^2
 \frac{\displaystyle\prod_{\nu=h+t+1}^{h+t+r}b_\nu
       \prod_{\nu=t+1}^{t+r}b_\nu}
      {\displaystyle\prod_{\nu=h+t+1}^{h+t+r+1}s_\nu(-\gamma)
       \prod_{\nu=t+1}^{t+r+1}s_\nu(\gamma)}.
\]
Here $\rho_r=t+r+1$. The prefix interval
$\{\rho_\nu:0\le \nu<r\}=\{t+1,\ldots,t+r\}$ shows that both products
are the corresponding factor of $\widetilde{\mathsf P}_r$ divided by
$s_{h+\rho_r}(-\gamma)s_{\rho_r}(\gamma)$. Consequently,
\begin{equation}\label{eq:Qtelescoping}
 (U_{-\gamma})_{ir}(U_\gamma)_{jr}-s_0(\gamma)^2(V_{-\gamma})_{ir}(V_\gamma)_{jr}
 =\frac{(-1)^{h+1}s_0(\gamma)^2\widetilde{\mathsf P}_r}{s_{h+\rho_r}(-\gamma)s_{\rho_r}(\gamma)}
 =\frac{(-1)^{h+1}s_0(\gamma)}{s_h(-\gamma)}\bigl(\widetilde{\mathsf P}_r-\widetilde{\mathsf P}_{r+1}\bigr),
\end{equation}
by \eqref{eq:oppositephasesine} with $a=h+\rho_r$, $a'=\rho_r$. Summation
gives $(Q_\gamma)_{ij}=(-1)^{h+1}s_0(\gamma)/s_h(-\gamma)$, which is $1$
for $h=0$. Transposing the definition of $Q_\gamma$ and using
$s_0(-\gamma)^2=s_0(\gamma)^2$ gives $Q_\gamma^T=Q_{-\gamma}$, which supplies
the upper triangle; equivalently, the formula holds for negative $h$
because $s_h(-\gamma)=-s_{-h}(\gamma)$.

\smallskip
\noindent\emph{Removing the temporary assumption $s_0(\gamma)\ne0$.}
When $s_\nu(\gamma)s_\nu(-\gamma)\ne0$ for $1\le \nu<n$, all entries of $U_{\pm\gamma}$,
$V_{\pm\gamma}$ (hence of $C_\gamma$ and $Q_\gamma$) are continuous in the
phase on the open set where this hypothesis holds, and so are the
right-hand sides of \eqref{eq:bilinearentries}, whose denominators involve
only $s_\nu(\pm\gamma)$ with $1\le \nu<n$. The phases with $s_0(\gamma)=0$ are
isolated, so the identities extend to them by continuity. This proves both formulas in \eqref{eq:bilinearentries}.

Taking $\gamma=\theta$ in the first formula and $\gamma=-\theta$ in the second,
and using $s_0(-\theta)=-s_0(\theta)$, gives \eqref{eq:six-C-entries}
and \eqref{eq:six-Q-entries}.
\end{proof}
\raggedbottom
\begingroup
\linespread{0.97}\selectfont
\providecommand{\bysame}{\leavevmode\hbox to3em{\hrulefill}\thinspace}
\providecommand{\MR}{\relax\ifhmode\unskip\space\fi MR }
\providecommand{\MRhref}[2]{%
  \href{http://www.ams.org/mathscinet-getitem?mr=#1}{#2}
}
\providecommand{\href}[2]{#2}

\endgroup
\end{document}